\documentclass[12pt,a4paper,twoside,final,notitlepage, reqno]{article}

\usepackage[english]{babel}
\usepackage[latin1]{inputenc}
\usepackage[a4paper,left=2.5cm,right=2.5cm]{geometry}
\usepackage{amssymb}
\usepackage{amsthm} 
\usepackage{amsmath}
\usepackage{amscd}
\usepackage{geometry}
\usepackage{graphics,graphicx}
\usepackage{epstopdf}
\usepackage[usenames, dvipsnames]{color}
\usepackage{hyperref}
\usepackage[textsize=small]{todonotes}

\usepackage[shortlabels]{enumitem}
\usepackage{subcaption}
\usepackage{siunitx}

\graphicspath{{./figures/}}

\theoremstyle{definition}

\theoremstyle{definition}

\theoremstyle{remark}
\newtheorem{remark}{Remark}

\providecommand{\keywords}[1]  {\textbf{Keywords:} #1}
\providecommand{\subjclass}[1] {\textbf{Subject classification:} #1}
\providecommand{\acknow}[1] {\textbf{Acknowledgments.} #1}

\renewcommand{\phi}{\varphi}

\usepackage{authblk}

\title{
\vspace{-3cm}
\bf{\Large{
  Analysis of time-stepping methods 
  \\[5pt]
  for the monodomain model
}}
}

\author[1]{Thomas Roy \thanks{thomas.roy@maths.ox.ac.uk}}
\author[2]{Yves Bourgault\thanks{ybourg@uottawa.ca}}
\author[3]{Charles Pierre \thanks{charles.pierre@univ-pau.fr}}

\affil[1]{
  Mathematical Institute, University of Oxford, United Kingdom.
}
\affil[2]{
  Department of Mathematics and Statistics, University of Ottawa,
  Ontario, Canada.
}
\affil[3]{
  Laboratoire  de Math\'ematiques et de leurs Applications,
  UMR CNRS 5142, \protect \\
  Universit\'e de Pau et des Pays de l'Adour, France.
  }

\usepackage{fancyhdr}

\begin{document} 
\date{25 May, 2020}

\maketitle

\begin{abstract}
To a large extent, the stiffness of the bidomain and monodomain models depends on the choice of the ionic model, which varies in terms of complexity and
realism. In this paper, we compare and analyze a variety of time-stepping methods: explicit or semi-implicit, operator splitting, exponential, and deferred correction methods. We compare these methods for solving the bidomain model coupled with three ionic models of varying complexity and stiffness: the phenomenological Mitchell-Schaeffer model, the more realistic Beeler-Reuter model, and the stiff and very complex ten Tuscher-Noble-Noble-Panfilov (TNNP) model. For each method, we derive absolute stability criteria of the spatially discretized monodomain model and verify that the theoretical critical time-steps obtained closely match the ones in numerical experiments. We also verify that the numerical methods achieve an optimal order of convergence on the model variables and derived quantities (such as speed of the wave, depolarization time), and this in spite of the local non-differentiability of some of the ionic models. The efficiency of the different methods is also considered by comparing computational times for similar accuracy. Conclusions are drawn on the methods to be used to solve the monodomain model based on the model stiffness and complexity, measured respectively by the eigenvalues of the model's Jacobian and the number of variables, and based on strict stability and accuracy criteria.
\end{abstract}

\noindent
\\[3pt]
\\[3pt]
\acknow{
The authors would like to thank the University of Ottawa for graduate scholarships to the first author, and the Natural Science and Engineering Research Council (NSERC) of Canada for a research grant to the second author. Research exchanges between France and Canada were funded by a grant from the Agence Nationale de la Recherche  of France
(ANR project HR-CEM no.\ 13-MONU-0004-01).
}

\section*{Introduction}
The modelling of the electrical activity of the heart offers an interesting perspective on the understanding of cardiac pathologies and its treatments. This subject has great potential in biomedical sciences, as experiments on living hearts require considerable resources and provide only a partial picture of the electrical activity of the heart. For instance, realistically simulating the behaviour of the heart reduces the necessity for these kinds of experiments. Considering that heart diseases are a leading cause of death in Western countries \cite{statcan}, modelling in cardiac electrophysiology has the potential to significantly impact society. Also, this area has achieved great strides in the past few years, because of improved models, calculation methods and computer power ~\cite{clayton2011,cain2011,winslow2000}.  

The mono- and bi- domain models allow for an adequate modelling of the electrical  activity in the heart tissue ~\cite{KeeSne2004,Sundnes2006}. These models form a complex system of partial  differential  equations (PDE)  that  are coupled with a system of ordinary differential equations (ODE) describing the ionic activity in cardiac cells. Solving these models numerically requires large computational time and computations are limited to few heartbeats on ventricles and atria.

The ODE systems describing ionic activity vary in terms of complexity depending on the physiological processes accounted for at the cell level. These ionic models, coupled with the mono- or bi- domain model, require discretization in space and time. This article focuses on various types of time-stepping methods: explicit, semi-implicit, operator splitting, exponential, and deferred correction methods. The monodomain model is a reaction-diffusion equation and solving the diffusion part of the model explicitly results in severe stability constraints on the time-step. Therefore, the use of semi-implicit or implicit methods is viewed as crucial by most authors \cite{Ethier2008,torabi2014,keener1998}. To achieve this, operator splitting methods have been considered by some authors \cite{Sundnes2006,schroll2007,trangenstein2004,hanslien2011}. Alternatively, an exponential method called the Rush-Larsen (RL) method was developed specifically to solve the ODEs resulting from the ionic activity of excitable cells, providing very stable numerical solutions \cite{RL1978}. This method is very popular and several other variations have been considered, including high-order RL methods ~\cite{PerVen2009,Marsh2013,sundnes2009}.

The different time-stepping methods have been thoroughly studied for non-spatial ionic models \cite{SpitDean2010,spiteri2008}, i.e. only at the cell level, and relatively little at the level of the heart tissue for the mono- or bi-domain models. To date, no research in cardiac electrophysiology has compared a large number of time-stepping methods thoroughly for spatial models. In addition, first-order numerical methods such as Euler's method are still widely used in this field of research due to their easy implementation \cite{puwal2007}. We seek to show that these first-order methods are very inefficient compared to high-order methods. Furthermore, we will consider a general ionic model for our theoretical stability analysis, allowing other researchers to use these results to determine critical time-steps for stability for the model of their choice. Earlier studies are specific to a given ionic model, which is limiting in a research area where models are constantly evolving. Thus, this research will provide a summary of the various methods that had not yet been compared, both for stability and accuracy.

The goal of this article is to find optimal time-stepping methods for ionic models of increasing stiffness, where optimality is judged based on the cost to reach a given accuracy. In Section \ref{chap:models}, we will start by choosing a number of ionic models that cover the range of complexity, stiffness and realism for these types of models.
In Section \ref{chap:methods}, we will list a number of viable numerical methods for the monodomain model. The different methods will be studied with theoretical analysis and numerical tests in the later section. 
In Section \ref{chap:stab}, different stability analysis will be conducted for each numerical method. We will establish a technique to study the stability of the numerical methods for the monodomain model in an ODE setting. To ease the comparison from one model to another, a general ionic model will be used. The theoretical time-steps will then be compared with the numerically observed critical time-step of the different methods. 
Then, in Section \ref{chap:convergence}, a convergence test will be performed to check whether the methods of high order exhibit the correct rate of convergence when used to solve the more complex ionic models. The accuracy of the methods studied will then be compared with respect to their accuracy, relatively to the size of the time-step and the computational time.

\section{Models in cardiac electrophysiology}\label{chap:models}

\subsection{Modelling the ionic activity in cells}\label{sec:cells}
There exist many different models used to describe the ionic activity in cardiac cells, varying in terms of complexity and stiffness. Here, we choose to study three models that cover the whole spectrum.

The Mitchell-Schaeffer (MS) model is a phenomenological two-variable model proposed in \cite{MS2003}, which can be easily understood analytically and is very efficient for numerical simulations. In spite of its very low complexity, this model still accurately represents the main characteristics of the cardiac action potential.

One of the physiological ionic models used in this paper is the Beeler-Reuter (BR) model \cite{BR1977}. This model has six gating variables and one concentration, besides the transmembrane potential. The source term in the potential equation includes four different ionic currents.

A more complex ionic model used in this article is the ten Tuscher-Noble-Noble-Panfilov (TNNP) model from \cite{TNNP2004}. There exist different versions of this model for different types of cells, and here we consider the epicardial cells. This version is the stiffest of the original TNNP models. This model has seventeen variables, including twelve gating variables and four concentrations. It also has 15 different ionic currents. 

Physiological models can be written in the following form, where $u$ is the transmembrane potential, $\mathbf{v}$ is the vector of gating variables, $\mathbf{X}$ is the vector of ionic concentrations, and $t$ is the time:
\begin{align}
\dfrac{d{u}}{d{t}} &=I(u,\mathbf{v},\mathbf{X},t), \label{eq:physmod1}\\
\dfrac{d{\mathbf{v}}}{d{t}} &= \mathbf{f}(u,\mathbf{v}),\label{eq:physmod2}\\
\dfrac{d{\mathbf{X}}}{d{t}} &= \mathbf{g}(u,\mathbf{v},\mathbf{X}) \label{eq:physmod3}.
\end{align}
The total current $I$ acts as a source/sink term and is defined as 
\begin{equation}\label{eq:sourcesink}
I(u,\mathbf{v},\mathbf{X},t)=\dfrac{1}{C_m} (I_\mathrm{app}(t)-I_\mathrm{ion}(u,\mathbf{v},\mathbf{X})),
\end{equation}
 where $C_m$ is the membrane capacity, $I_\mathrm{app}$ is the applied stimulation current, and $I_\mathrm{ion}$ is the sum of all ionic currents across the cell membrane. In general, the r.h.s. $f_i$'s in (\ref{eq:physmod2}) mimic Hodgkin-Huxley type equations \cite{HodgHux1952}:
\begin{equation}\label{eq:hodhux}
f_i(u,v_i)=\dfrac{v_{i,\infty}(u) -v_i}{\tau_{v_i}(u)},
\end{equation}
where $v_{i,\infty}$ is a steady state value for $v_i$ and $\tau_{v_i}$ is a time constant, both specific to the ionic model and the gating variable. The function $\mathbf{g}$ in (\ref{eq:physmod3}) is specific to the ionic model. 

\subsection{Modelling the Heart Tissue}

After modelling the electrical activity at the cell level, we now extend our models to the tissue level. For our test cases, we consider the monodomain model with constant extra- and intra-cellular conductivities,  
coupled with a general ionic cell model. For the monodomain model, we assume equal anisotropy ratio of the conductivities.
The monodomain model reads as:
\begin{align}
\label{eq:nondim01}
\frac{\partial u}{\partial t} &= I(u,\mathbf{v},\mathbf{X},t,x) + \mathrm{div} (\sigma \nabla u),
\\
\label{eq:nondim02}
\frac{\partial \mathbf{v}}{\partial t}&=\mathbf{f}(u,\mathbf{v}),
\\
\label{eq:nondim03}
\frac{\partial \mathbf{X}}{\partial t}&=\mathbf{g}(u,\mathbf{v},\mathbf{X}),
\end{align}
where $\sigma=\frac{\lambda}{1+\lambda}\sigma_i /\chi C_m$, $\lambda$ is the equal anisotropy ratio, $\sigma_i$ is the intra-cellular conductivity tensor and $\chi$ is the cellular membrane area per unit volume of cardiac tissue. The applied stimulation current $I_\mathrm{app}$ in \eqref{eq:sourcesink} can now vary through space, denoted by $x\in \Omega$, the spatial domain. We denote by $\mathbf{f}$ and $\mathbf{g}$ the vectors of $f_i$'s and $g_i$'s, respectively. In some cases, it is simpler to group in $\mathbf{v}$ the vector of $v_i$'s and $X_i$'s, and in $\mathbf{f}$ the vector of $f_i$'s and $g_i$'s.

We impose homogeneous Neumann boundary conditions, i.e. $\nabla u\cdot n =0$ on $\partial \Omega$. This amounts to having no current leaking from the heart to the surrounding tissues. The initial condition is taken as the resting state of the chosen ionic model.

\section{Numerical methods}\label{chap:methods}
We consider a variety of explicit and semi-implicit methods of first to third-order, using a constant time-step $\Delta t$, for the time discretization. When solving reaction-diffusion equations as the monodomain model, some authors suggest using semi-implicit methods with implicit diffusion and explicit reaction terms~\cite{Ethier2008,PerVen2009,Sundnes2006}. Implicit methods are usually too expensive for many ionic models and explicit methods have additional stability restrictions. 

Let us introduce the notation for the fully discretized version of the problem  (\ref{eq:nondim01})-(\ref{eq:nondim03}). Given a spatial mesh with nodes $x_i$, $i=0,1,\dots, N$ and a temporal mesh with equally spaced temporal nodes $t_n= n\Delta t$, $n\geq 0$, we denote by
\begin{align}U^n &=[U_0^n,\ldots ,U_{N}^n]^\top \simeq  [u(x_0,t_n),\ldots ,u(x_{N},t_n)]^\top ,\\
 V^n &=[V_{1,0}^n,\ldots ,V_{p,0}^n,\ldots ,V_{p,N}^n]^\top \simeq  [v_1(x_0,t_n),\ldots ,v_p(x_0,t_n),\ldots ,v_p(x_{N},t_n)]^\top , \\
X^n &=[X_{1,0}^n,\ldots ,X_{q,0}^n,\ldots ,X_{q,N}^n]^\top \simeq  [X_1(x_0,t_n),\ldots ,X_q(x_0,t_n),\ldots ,X_q(x_{N},t_n)]^\top ,\end{align}
the approximate values for $u$, $\mathbf{v}$ and $\mathbf{X}$, respectively. 

We write the equations for the finite difference and finite element methods. The term $AU^n$ is the discrete version of $\mathrm{div} (\sigma \nabla u)$. For the finite difference method, $A$ is a discretized Laplacian, and for the finite element method, $A=M^{-1} S$, where $M$ is the mass matrix and $S$ is the stiffness matrix.

For most methods, the gating variables and concentrations are treated the same way. Therefore, if the following schemes do not have $X$ and $G$, it is implied they are included in $V$ and $F$, respectively.

\subsection{First-order methods}
\begin{enumerate}[(i)]
\item Forward Euler (FE):
\begin{equation}\label{eq:FEscheme}
\begin{gathered}
\frac{U^{n+1}-U^n}{\Delta t}=I(U^n,V^n,t_n,x)+AU^n,\\
\frac{V^{n+1}-V^n}{\Delta t}=F(U^n,V^n).
\end{gathered}
\end{equation}
\item Forward-Backward Euler (FBE):
\begin{equation}\label{eq:FBEscheme}
\begin{gathered}
\frac{U^{n+1}-U^n}{\Delta t}=I(U^n,V^n,t_n,x)+AU^{n+1},\\
\frac{V^{n+1}-V^n}{\Delta t}=F(U^n,V^n).
\end{gathered}
\end{equation}
\item Rush-Larsen with Forward-Backward Euler (RL-FBE) ~\cite{RL1978,PerVen2009}:
\begin{equation}\label{eq:RLscheme} 
\begin{gathered}
\frac{U^{n+1}-U^n}{\Delta t}=I(U^n,V^n,X^n,t_n,x)+AU^{n+1},\\
\frac{V_i^{n+1}-V_i^n}{\Delta t}=\Phi (a_i^n\Delta t)(a_i^n V_i^n+b_i^n),\quad  i=1,\ldots ,p,\\
\frac{X^{n+1}-X^n}{\Delta t}=G(U^n,V^n,X^n),
\end{gathered}
\end{equation}
where $\Phi, a_i^n$ and $b_i^n$ are given below. This method corresponds to the case $c_{-1}=0$, $c_0=1$, $c_1=0$ in the general RL methods described at the end of this section.
\end{enumerate}
\subsection{Second-Order Methods}
\begin{enumerate}[(i)]
\item Second-order semi-implicit backward differentiation (SBDF2) ~\cite{Ascher1995,Ethier2008}:
\begin{equation}\label{eq:SBDFscheme}
\begin{gathered}
\frac{\frac{3}{2}U^{n+1}-2U^n+\frac{1}{2}U^{n-1}}{\Delta t} = 2I(U^n,V^n,t_n,x) - I(U^{n-1},V^{n-1},t_{n-1},x)+AU^{n+1},\\
\frac{\frac{3}{2}V^{n+1}-2V^n+\frac{1}{2}V^{n-1}}{\Delta t} = 2F(U^n,V^n) - F(U^{n-1},V^{n-1}).
\end{gathered}
\end{equation}

\item Strang Splitting with Crank-Nicolson and Runge-Kutta 2 (CN-RK2) ~\cite{Strang1968,Sundnes2006}:\\
For this method, we do half a time-step of RK2 solely on the reaction part of the monodomain equations, followed by a step of CN on the diffusion part and another half-step of RK2 on the reaction part. Here we denote $Y^n=\begin{bmatrix} U^{n} \\ V^{n} 
\end{bmatrix}$, $Y^*=\begin{bmatrix} U^{*} \\ V^{*} 
\end{bmatrix}$ and $Y^{**}=\begin{bmatrix} U^{**} \\ V^{**} 
\end{bmatrix}$.\\
\textbf{Step 1:}
\begin{equation}\label{eq:cnrk2scheme1}
\frac{Y^*-Y^n}{\Delta t /2}=\begin{bmatrix} I\left(Y^n+\frac{\Delta t}{4}\begin{bmatrix} I(Y^n,t_n,x) \\ F(Y^n) \end{bmatrix},t_n+\frac{\Delta t}{4},x\right) \\ F\left(Y^n+\frac{\Delta t}{4}\begin{bmatrix} I(Y^n,t_n,x) \\ F(Y^n) \end{bmatrix}\right)\end{bmatrix}.
\end{equation}
\textbf{Step 2:}
\begin{equation}\label{eq:cnrk2scheme2}
\frac{U^{**}-U^*}{\Delta t}=\frac{1}{2}A(U^{**}+U^*).
\end{equation}
\textbf{Step 3:}
\begin{equation}\label{eq:cnrk2scheme3}
\frac{Y^{n+1}-Y^{**}}{\Delta t /2}=\begin{bmatrix} I\left(Y^{**}+\frac{\Delta t}{4}\begin{bmatrix} I(Y^{**},t_n+\Delta t/2,x) \\ F(Y^{**}) \end{bmatrix},t_n+\frac{3}{4}\Delta t,x\right) \\ F\left(Y^{**}+\frac{\Delta t}{4}\begin{bmatrix} I(Y^{**},t_n+\Delta t/2,x) \\ F(Y^{**}) \end{bmatrix}\right) 
\end{bmatrix}.
\end{equation}

\item Strang Splitting with Crank-Nicolson and Runge-Kutta 4 (CN-RK4):\\
This scheme is written as the CN-RK2 scheme, but using instead the fourth-order Runge-Kutta method (see \cite{Roy2015} for more details).

\item Second-order Rush-Larsen with Crank-Nicolson Adam-Bashforth (RL-CNAB) ~\cite{RL1978,PerVen2009}:
\begin{equation}\label{eq:RL2scheme}
\begin{gathered}
\frac{U^{n+1}-U^n}{\Delta t} = \frac{3}{2}I(U^n,V^n,X^n,t_n,x) - \frac{1}{2}I(U^{n-1},V^{n-1},X^{n-1},t_{n-1},x)+\frac{1}{2}A(U^{n+1}+U^n),\\
\frac{V_i^{n+1}-V_i^n}{\Delta t}=\Phi (a_i^{n+\frac{1}{2}}\Delta t)(a_i^{n+\frac{1}{2}} V_i^n+b_i^{n+\frac{1}{2}}),\quad  i=1,\ldots ,p,\\
\frac{X^{n+1}-X^n}{\Delta t}=\frac{3}{2}G(U^n,V^n,X^n)-\frac{1}{2}G(U^{n-1},V^{n-1},X^{n-1}),
\end{gathered}
\end{equation}
where $\Phi, a_i^n$ and $b_i^n$ are given below. This method corresponds to the case $c_{-1}=0$, $c_0=\frac{3}{2}$, $c_1=-\frac{1}{2}$ in the general RL methods described at the end of this section.
\end{enumerate}

\subsection{Third-Order Methods}
\begin{enumerate}[(i)]
\item Third-Order Deferred Correction (DC3) ~\cite{KreGus2002,MinevGuermond}:\\
Here we denote $Y^n = \begin{bmatrix} U^{n} \\ V^{n} \end{bmatrix}$ and  $Y_i^n = \begin{bmatrix} U_i^{n} \\ V_i^{n} \end{bmatrix}$ .  The initial values for the substeps are $Y_0^0 = Y_0$, $Y_1^0 = 0$,  and $Y_2^0 = 0$. This method proceeds with three substeps:
\begin{equation}\label{eq:DCn0}
\text{for } n\geq 0\quad \begin{cases}
				ml_0^{n+1} = I(Y_0^n,t_n,x), \quad nl_0^{n+1} = F(Y_0^n),\\
				\dfrac{U_0^{n+1} - U_0^n}{\Delta t} = ml_0^{n+1}+ AU_0^{n+1}, \quad \dfrac{V_0^{n+1} - V_0^n}{\Delta t} =nl_0^{n+1},	\\
				dY_0^{n+1} = (Y_0^{n+1} - Y_0^n)/\Delta t
			\end{cases}
\end{equation}

\begin{equation}\label{eq:DCn1}
\text{for } n\geq 1\quad \begin{cases}
				d^2U_0^{n+1} = (dU_0^{n+1} - dU_0^n)/ \Delta t, \quad ml_1^{n} = I(Y_0^n + \Delta t Y_1^{n-1},t_n,x),\\
				\dfrac{U_1^{n} - U_1^{n-1}}{\Delta t} = AU_1^{n} - \dfrac{1}{2}d^2U_0^{n+1} + \dfrac{ml_1^n -  ml_0^{n}}{\Delta t} ,	\\
				d^2V_0^{n+1} = (dV_0^{n+1} - dV_0^n)/ \Delta t, \quad nl_1^{n} = F(Y_0^n + \Delta t Y_1^{n-1}),\\
				\dfrac{V_1^{n} - V_1^{n-1}}{\Delta t} = -\dfrac{1}{2}d^2V_0^{n+1} + \dfrac{nl_1^n -  nl_0^{n}}{\Delta t} ,	\\
				dY_1^{n} = (Y_1^{n} - Y_1^{n-1})/\Delta t,
			\end{cases}
\end{equation}

\begin{equation}\label{eq:DCn2}
\text{for } n\geq 2\quad \begin{cases}
				d^2U_1^{n} = (dU_1^{n} - dU_1^{n-1})/ \Delta t, \quad	d^3U_0^{n+1} = (d^2U_0^{n+1} - d^2U_0^{n})/ \Delta t, \\
				 ml_2^{n-1} = I(Y_0^{n-1} + \Delta t Y_1^{n-1} + \Delta t^2 Y_2^{n-2},t_n,x),\\
				\dfrac{U_2^{n-1} - U_2^{n-2}}{\Delta t} = AU_2^{n-1} - \dfrac{1}{2}d^2U_1^{n} + \dfrac{1}{6}d^3U_0^{n+1} + \dfrac{ml_2^{n-1} -  ml_1^{n-1}}{\Delta t^2} ,\\
				d^2V_1^{n} = (dV_1^{n} - dV_1^{n-1})/ \Delta t \quad	d^3V_0^{n+1} = (d^2V_0^{n+1} - d^2V_0^{n})/ \Delta t, \\
				 nl_2^{n-1} = F(Y_0^{n-1} + \Delta t Y_1^{n-1} + \Delta t^2 Y_2^{n-2}),\\
				\dfrac{V_2^{n-1} - V_2^{n-2}}{\Delta t} = -\dfrac{1}{2}d^2V_1^{n} + \dfrac{1}{6}d^3V_0^{n+1} + \dfrac{nl_2^{n-1} -  nl_1^{n-1}}{\Delta t^2} ,\\
				Y^{n-1} = Y_0^{n-1} + \Delta t Y_1^{n-1} + \Delta t^2 Y_2^{n-1}.
			\end{cases}
\end{equation}
It is possible to obtain a second order method by solving the first two substeps and defining $Y^{n-1}= Y_0^{n-1} + \Delta t Y_1^{n-1}$. We only provide numerical results for the third-order DC method. 

\end{enumerate}

\subsection{The Rush-Larsen Method}

A very popular method for solving physiological models is the first-order scheme proposed by Rush and Larsen \cite{RL1978}. Outside of cardiac electrophysiology, this method is known as the explicit exponential Euler method. The Rush-Larsen (RL) method solves problems of the form $\frac{dy}{dt} = 	a(t)y + b(t)$.

Perego and Veneziani \cite{PerVen2009} introduced a way to increase the order of the scheme, by taking $a$ and $b$ at time $t_{n+\frac{1}{2}}$,
\begin{equation}\label{eq:RL2}
\begin{cases}
y^{n+1}=y^n+\Delta t\Phi (a^{n+\frac{1}{2}}\Delta t)(a^{n+\frac{1}{2}}y^n +b^{n+\frac{1}{2}}),\quad n=0,\ldots , N_t,\\
y(0)=y^0,
\end{cases}
\end{equation}
where
\begin{equation*}
\Phi (x)=	\begin{cases}
		\dfrac{e^x-1}{x},\quad&x\neq 0,\\
		1,\quad&x=0,
		\end{cases}
\end{equation*}
$ a^{n+\frac{1}{2}}$ and $b^{n+\frac{1}{2}}$ are approximations of $a(t_{n+\frac{1}{2}})$ and $b(t_{n+\frac{1}{2}})$ taken as 
\begin{equation}
\begin{aligned}
 &a^{n+\frac{1}{2}}=c_{-1}a^{n+1}+c_0 a^n +c_{1}a^{n-1}, \quad b^{n+\frac{1}{2}}=c_{-1}b^{n+1}+c_0 b^n +c_{1}b^{n-1}, \;n=1,\ldots ,N_t,\\
 &a^{\frac{1}{2}}=c_{-1}a^1+(c_0+c_{1})a^0, \quad b^{\frac{1}{2}}=c_{-1}b^1+(c_0+c_{1})b^0,
\end{aligned}
\end{equation}
where $c_{-1}$, $c_0$ and $c_1$ are coefficients to be determined. 

For the transmembrane potential $u$ and concentrations $\mathbf{X}$, we take $a=0$, which simply results in the Euler and Adams-Bashforth methods for the first and second order cases, respectively. In each case we also solve the diffusion implicitly, hence resulting in the FBE and CNAB methods for the transmembrane potential and concentrations.

For the gating variables, we write the function $\mathbf{f}$ in  (\ref{eq:hodhux}) as $f_i(u,\mathbf{v},\mathbf{X})=a_i(u)v_i +b_i(u)$. For the BR and TNNP models, we have:
\begin{equation}
a_i(u)=\dfrac{-1}{\tau_{v_i}(u)},\quad\text{and}\quad b_i(u)=\dfrac{v_{i,\infty}(u)}{\tau_{v_i}(u)}.\end{equation}

\section{Stability analysis}\label{chap:stab}
We now find stability conditions for the time-stepping methods introduced in the previous section. 

\subsection{Absolute Stability}
The stability analysis will be done using the concept of absolute stability for ODE solvers. We linearize our PDE and then consider the linear ODE system resulting from the linearized problem semi-discretized in space. Since we will carry von Neumann stability analysis, we will consider 1D reaction-diffusion equations on the whole real line discretized in space by centred finite difference schemes. We find the stability region of each method for the linearized problem.  This analysis will confirm that the stability conditions of semi-implicit methods applied to the monodomain model coincide with their explicit equivalent applied solely to the ionic models. The largest possible time-step $\Delta t$ for which $\lambda \Delta t$ is in the stability region for all $\lambda$ is called the {\em critical} time-step, where $\lambda$ are the eigenvalues of the linearized problem.

Considering (\ref{eq:nondim01})-(\ref{eq:nondim03}), but denoting by $\mathbf{v}$ the vector of $v_i$'s and $X_i$'s, and $\mathbf{f}$ the vector of $f_i$'s and $g_i$'s, we linearize around some state $(\tilde{u},\tilde{\mathbf{v}})$. This state is constant in space and time, and will be determined later. We obtain
\begin{equation}\label{eq:lin1}
\frac{\partial u}{\partial t} = \frac{\partial I}{\partial u}(\tilde{u},\tilde{\mathbf{v}})u + \frac{\partial I}{\partial \mathbf{v}}(\tilde{u},\tilde{\mathbf{v}})\mathbf{v} +\mathrm{div}(\sigma\nabla u),
\end{equation}
\begin{equation}\label{eq:lin2}
\frac{\partial \mathbf{v}}{\partial t}=\frac{\partial \mathbf{f}}{\partial u}(\tilde{u},\tilde{\mathbf{v}})u+\frac{\partial \mathbf{f}}{\partial \mathbf{v}}(\tilde{u},\tilde{\mathbf{v}})\mathbf{v}.
\end{equation}
For the sake of the stability analysis, we consider uniformly spaced grid points $x_{j}=jh$ for each $j\in\mathbb{Z}$, $h>0$, the space step. We denote \begin{align}
U &=[\ldots ,U_{-1},U_0,U_1,\ldots ]^\top \simeq  [\ldots ,u(x_{-1},t),u(x_0,t),u(x_1,t),\ldots ]^\top,\\
V &=[\ldots ,V_{-1},V_{0},V_{1},\ldots ]^\top \simeq [\ldots ,\mathbf{v}(x_{-1},t),\mathbf{v}(x_0,t),\mathbf{v}(x_1,t),\ldots ]^\top.
\end{align}
Discretizing (\ref{eq:lin1}) and (\ref{eq:lin2}) we obtain
\begin{equation}\label{eq:discret}
\dfrac{d}{dt}\begin{bmatrix} U_j \\ V_j \end{bmatrix}=
 \begin{bmatrix} \dfrac{\partial I}{\partial u}(\tilde{u},\tilde{\mathbf{v}}) & \dfrac{\partial I}{\partial \mathbf{v}}(\tilde{u},\tilde{\mathbf{v}}) 
\\[0.6em]
\dfrac{\partial \mathbf{f}}{\partial u}(\tilde{u},\tilde{\mathbf{v}}) & \dfrac{\partial \mathbf{f}}{\partial \mathbf{v}}(\tilde{u},\tilde{\mathbf{v}})\end{bmatrix}
\begin{bmatrix} U_j \\ V_j \end{bmatrix}+\sigma\begin{bmatrix} \dfrac{U_{j+1}-2U_j+U_{j-1}}{h^2} \\[0.5em] 0\end{bmatrix},
\end{equation}
for each $j\in\mathbb{Z}$. As is usually done for von Neumann stability analysis,
we set $U_j (t)= U_\omega (t)e^{i\omega jh}$ and $V_j (t)= V_\omega(t)e^{i\omega jh}$, for all $\omega \in [0,\dfrac{2\pi}{h}]$. Then
\begin{equation}\label{eq:fourier}
\dfrac{d}{dt}\begin{bmatrix} U_\omega \\ V_\omega \end{bmatrix}=
(J+A_\omega)
\begin{bmatrix} U_\omega \\ V_\omega \end{bmatrix},
\end{equation}
where $J$ is the Jacobian matrix evaluated at $(\tilde{u},\tilde{\mathbf{v}})$ as in (\ref{eq:discret}) and $A_\omega$ is a diffusion matrix given by
\begin{equation} \label{eq:Aomega}
A_\omega = \frac{2\sigma}{h^2}\begin{bmatrix} \cos(\omega h)-1 & 0 \\ 0 & 0 \end{bmatrix}. 
\end{equation}
Equation (\ref{eq:fourier}) is discretized with the time-stepping methods of Section~\ref{chap:methods}.
A stability function $R=R(\Delta t, h, \omega))$ is then defined for each method. The condition for stability is that the spectral radius $\rho(R)\le 1$, assuming that all the eigenmodes of $J+A_\omega$ satisfy ${\rm Re}(\lambda)<0$.

\subsubsection{Forward Euler}

We now discretize in time and apply the FE scheme to (\ref{eq:fourier}). We obtain
\begin{equation}
\begin{bmatrix} U_\omega^{n+1} \\ V_\omega^{n+1} \end{bmatrix}=(I_d +
\Delta t (J+A_\omega))\begin{bmatrix} U_\omega^{n} \\ V_\omega^{n}\end{bmatrix},
\end{equation}
where $I_d$ is the identity matrix.

The stability condition is $|1+\Delta t\lambda|\le 1$ for all $\lambda$ eigenvalues of $(J+A_\omega)$. Computing numerically the spectrum of $(J+A_\omega)$, we noticed that this condition is the most stringent for the most negative eigenvalue $\lambda_{\mathrm{min}}$ of $J+A_\omega$. The condition for stability is then $\Delta t \leq -2/\lambda_{\mathrm{min}}$. Because this eigenvalue depends on $\omega$, we will consider the wave number $\omega$ that makes this inequality most restrictive (in the sense that it makes $-2/\lambda_{\mathrm{min}}$ the smallest), i.e. $\omega=\pi/h$, a fact that can be verified numerically. Therefore our critical time-step for the stability of the FE method is $\Delta t^*_{\mathrm{theo}}=-2/\lambda_{\mathrm{min}}$, where $\lambda_{\mathrm{min}}$ is the minimum eigenvalue of 
\begin{equation}
\begin{bmatrix} \dfrac{\partial I}{\partial u}(\tilde{u},\tilde{\mathbf{\mathbf{v}}}) -\dfrac{4\sigma}{h^2}& \dfrac{\partial I}{\partial \mathbf{v}}(\tilde{u},\tilde{\mathbf{v}}) 
\\[0.6em] \dfrac{\partial \mathbf{f}}{\partial u}(\tilde{u},\tilde{\mathbf{v}}) & \dfrac{\partial \mathbf{f}}{\partial \mathbf{v}}(\tilde{u},\tilde{\mathbf{v}})\end{bmatrix}.
\end{equation}
\\
To determine $\lambda_{\mathrm{min}}$, we calculate numerically the solution of the problem on a fine grid. We then evaluate the matrix $J+A_\omega$ at each node of our domain, and for each of these matrices, we calculate their eigenvalues. We choose the most negative eigenvalue as our $\lambda_{\mathrm{min}}$. 
The constant state $(\tilde{u},\tilde{\mathbf{v}})$ is then chosen as the numerical solution evaluated at this node. The most negative eigenvalue typically appears during the resting state of the solution, but for some models, the stimulation current can sometimes make this eigenvalue more negative. The theoretical critical time-steps are shown in Table \ref{tab:crittheo} for the BR model, Table \ref{tab:crittheoMS} for the MS model and in Table \ref{tab:FEtnnp} for the TNNP model. For small values of $h$, the critical time-step is proportional to $h^2$.

\subsubsection{Forward-Backward Euler}

We define the stability function
$R(\Delta t, h, \omega)=(I-\Delta t A_\omega)^{-1}(I_d+\Delta t J).$
One can easily see that the matrix $(I_d-\Delta t A_\omega)^{-1}$ is diagonal with value $\left(1-\Delta t \sigma \dfrac{2 \cos \omega h -2}{h^2}\right)^{-1}$ in the first position and 1 on the rest of the diagonal, implying that the spectral radius $\rho((I_d-\Delta t A_\omega)^{-1}) = 1$.  This yields  $\rho(R(\Delta t, h, \omega)) \le \rho(I_d+\Delta t J)$, which inequality is sharp when $\omega=0$. This removes the dependence on $h$ in the stability condition. 

The condition for stability then becomes $|1+\Delta t\lambda|\le 1$ for all $\lambda$ eigenvalues of $J$. Therefore, similarly to the FE method, our critical time-step for FBE is $\Delta t^*_{\mathrm{theo}}=-2/\lambda_{\mathrm{min}}$, where $\lambda_{\mathrm{min}}$ is the most negative eigenvalue of $J$. As explained for the previous method, we calculate $J$ at each node of a fine grid and we choose its most negative eigenvalue as our $\lambda_{\mathrm{min}}$. The constant state $(\tilde{u},\tilde{\mathbf{v}})$ is then chosen as the numerical solution evaluated at that node. The critical time-steps for the BR, MS and TNNP models are shown in Tables \ref{tab:crittheo}, \ref{tab:crittheoMS} and \ref{tab:crittheoTNNP}, respectively. 

On coarse meshes, FE and FBE are expected to be stable for the same time-step $\Delta t$, but as the grid is refined, FE's stability deteriorates while FBE's remains unchanged. In fact, the diffusion term taken implicitly in the FBE method makes the critical time-step independent of the grid size $h$. 

\subsubsection{Strang Splitting}

We use the technique presented above to look at the stability of the Strang splitting scheme. Let us start with the CN-RK2 method.\\
\textbf{Step 1:}
\begin{equation}\label{eq:cnrk1}
\begin{bmatrix} U_\omega^{*} \\ \mathbf{V}_\omega^{*} \end{bmatrix} =\left(I_d+\frac{\Delta t}{2}J+\frac{(\Delta t/2)^2}{2}J^2 \right)\begin{bmatrix} U_\omega^{n} \\ \mathbf{V}_\omega^{n}\end{bmatrix}.
\end{equation}
\textbf{Step 2:}
\begin{equation}\label{eq:cnrk2}
\begin{bmatrix} U_\omega^{**} \\ \mathbf{V}_\omega^{**} \end{bmatrix}=
\left(I_d-\frac{1}{2}\Delta t A_\omega\right)^{-1}\left(\frac{1}{2}\Delta t A_\omega+I_d\right)
\begin{bmatrix} U_\omega^{*} \\ \mathbf{V}_\omega^{*}\end{bmatrix} 
 .
\end{equation}
\textbf{Step 3:}

\begin{equation}\label{eq:cnrk3}
\begin{bmatrix} U_\omega^{n+1} \\ \mathbf{V}_\omega^{n+1} \end{bmatrix} =\left(I_d+\frac{\Delta t}{2}J+\frac{(\Delta t/2)^2}{2}J^2 \right)\begin{bmatrix} U_\omega^{**} \\ \mathbf{V}_\omega^{**}\end{bmatrix}
.
\end{equation}
Combining (\ref{eq:cnrk1}), (\ref{eq:cnrk2}) and (\ref{eq:cnrk3}) we obtain
\begin{multline}
\begin{bmatrix} U_\omega^{n+1} \\ \mathbf{V}_\omega^{n+1} \end{bmatrix} =\left(I_d+\frac{\Delta t}{2}J+\frac{(\Delta t/2)^2}{2}J^2\right)\left(I_d-\frac{1}{2}\Delta t A_\omega\right)^{-1}\left(\frac{1}{2}\Delta t A_\omega+I_d\right)\\\left(I_d+\frac{\Delta t}{2}J+\frac{(\Delta t/2)^2}{2}J^2 \right)\begin{bmatrix} U_\omega^{n} \\ \mathbf{V}_\omega^{n}\end{bmatrix}
.
\end{multline}
For stability we need $\rho(R(\Delta t, h, \omega))\leq 1$. Noticing that
\begin{equation}
 \rho\left(\left(I_d-\frac{1}{2}\Delta t A_\omega\right)^{-1}\left(\frac{1}{2}\Delta t A_\omega+I_d\right)\right)\le 1,
\end{equation}
the condition for stability becomes
\begin{equation}\label{eq:strangstab}
\rho \left(I_d+\frac{\Delta t}{2}J+\frac{(\Delta t/2)^2}{2}J^2\right)^2\leq 1,
\end{equation}
independently of $h$.

\begin{figure}
\centering
\begin{subfigure}[b]{0.49\textwidth}
\includegraphics[width=.99\linewidth]{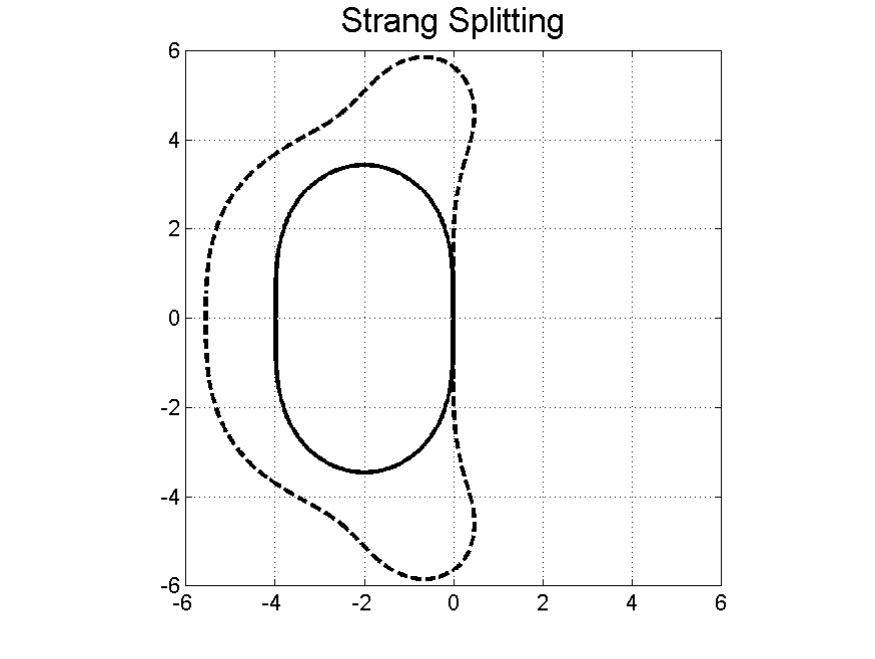}
 
  \caption{CN-RK4 (dotted), CN-RK2 (solid)} \label{fig:stabstrang}
\end{subfigure}
\begin{subfigure}[b]{0.49\textwidth}
\includegraphics[width=.99\linewidth]{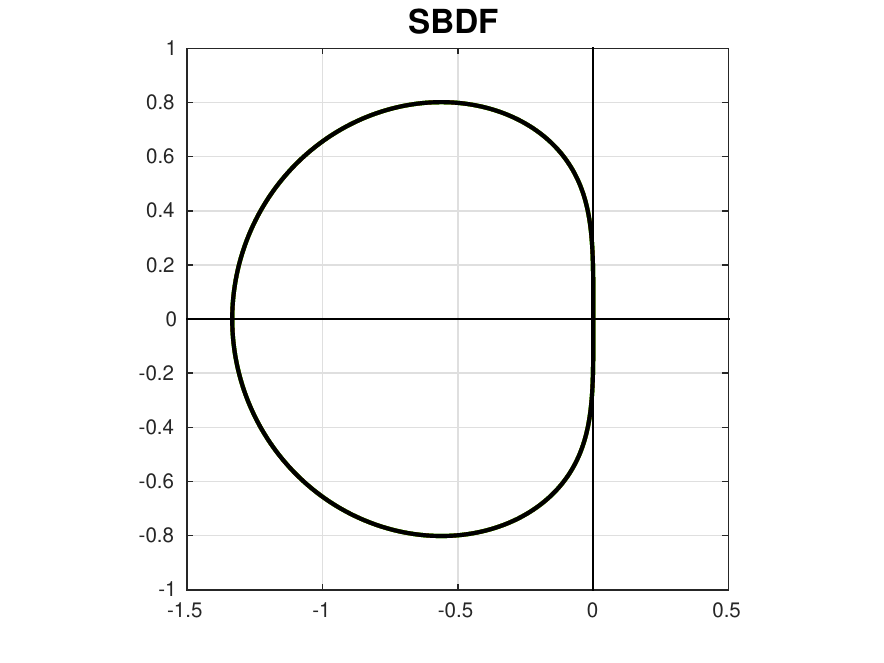}
 
  \caption{SBDF2 method} \label{fig:stabsbdf}
\end{subfigure}
\caption{Stability Regions}
\end{figure}

If $J$ is diagonalizable, (\ref{eq:strangstab}) is equivalent to
\begin{equation}\label{eq:cnrkstab2}
 \left|1+\Delta t \lambda_i/2+\dfrac{(\Delta t \lambda_i/2)^2}{2}\right|\leq 1,
\end{equation}
for all eigenvalue $\lambda_i$ of $J$. As seen with the previous methods, and with the contour of the stability region shown in Figure \ref{fig:stabstrang}, the most restrictive eigenvalue will be the most negative. Using this most negative eigenvalue we obtain the critical time-steps shown in Tables \ref{tab:crittheo}, \ref{tab:crittheoMS} and \ref{tab:crittheoTNNP}.

With a similar derivation, a condition for stability of the CN-RK4 method is obtained:
\begin{equation}\label{eq:splittstab2}
 \left|1+\Delta t \lambda_i/2+\dfrac{(\Delta t \lambda_i/2)^2}{2}+\dfrac{(\Delta t \lambda_i/2)^3}{6}+\dfrac{(\Delta t \lambda_i/2)^4}{24}\right|\leq 1,
\end{equation}
for all eigenvalue $\lambda_i$ of $J$. Again, we solve for the stability contour and get the critical time-steps $\Delta t^*_{\mathrm{theo}}$ in Tables \ref{tab:crittheo}, \ref{tab:crittheoMS} and \ref{tab:crittheoTNNP}. The stability region of CN-RK4 is shown in Figure \ref{fig:stabstrang} and contains the stability region of CN-RK2.

\subsubsection{Second-Order Semi-Implicit Backward Differentiation}

Following the approach used previously, we write the SBDF2 scheme as:
\begin{equation}\label{eq:SBDF1}
\dfrac{1}{\Delta t}\left(\frac{3}{2}\begin{bmatrix} U_\omega^{n+1} \\ V_\omega^{n+1} \end{bmatrix}-2\begin{bmatrix} U_\omega^{n} \\ V_\omega^{n} \end{bmatrix}+\frac{1}{2}\begin{bmatrix} U_\omega^{n-1} \\ V_\omega^{n-1} \end{bmatrix}\right)=
J\left(
2\begin{bmatrix} U_\omega^{n} \\ V_\omega^{n}\end{bmatrix}-\begin{bmatrix} U_\omega^{n-1} \\ V_\omega^{n-1} \end{bmatrix}\right)+A_\omega\begin{bmatrix} U_\omega^{n+1} \\ V_\omega^{n+1}\end{bmatrix}.
\end{equation}
For simplicity, we denote $Y^j=\begin{bmatrix} U_\omega^{j} \\ V_\omega^{j} \end{bmatrix}$, for all $j\in\mathbb{Z}$. Equation (\ref{eq:SBDF1}) becomes

\begin{equation}\label{eq:SBDF2}
(\frac{3}{2}I_d-A_\omega)Y^{n+1}-(2I_d+2\Delta t J)Y^n+(\frac{1}{2}I_d+\Delta t J)Y^{n-1}=0.
\end{equation}
As done for the previous semi-implicit methods, we take $\omega = 0$, which leads to a stability condition independent of $h$. See the remark below for a proof that this value is the most restrictive for the stability of this method. Equation (\ref{eq:SBDF2}) becomes
\begin{equation}\label{eq:SBDF3}
\frac{3}{2}Y^{n+1}-(2I_d +2\Delta t J)Y^n+(\frac{1}{2}I_d +\Delta t J)Y^{n-1}=0.
\end{equation}
Equation (\ref{eq:SBDF3}) is solved using the Lagrange's method found in \cite{Hairer1993}. We set $Y^j=\zeta^j W$, where $W$ is an eigenvector of $J$ with eigenvalue $\lambda$. We divide by $\zeta^{n-1}$ and obtain the following equation
\begin{equation}\label{eq:SBDF4}
\frac{3}{2}\zeta^2 W-(2I_d +2\Delta t J)\zeta W+(\frac{1}{2}I_d+\Delta t J)W=0,
\end{equation}
which implies
\begin{equation}\label{eq:SBDF5}
\frac{3}{2}\zeta^2-(2+2\Delta t \lambda)\zeta +(\frac{1}{2}+\Delta t \lambda)=0.
\end{equation}
Equation (\ref{eq:SBDF3}) has stable solutions \textit{iff} for any eigenvalue $\lambda$, all simple roots $\zeta (\lambda \Delta t)$ of (\ref{eq:SBDF5}) satisfy $|\zeta (\lambda \Delta t)|\leq 1$, and additional multiple roots must satisfy $|\zeta(\lambda \Delta t)|<1$ \cite{Hairer1996}. The stability region of the method is defined as 
\begin{equation}
S=\left\{\mu \in \mathbb{C}\; ; \begin{split} & \text{ all simple roots $\zeta(\mu)$ of (\ref{eq:SBDF5}) satisfy $|\zeta(\mu)|\leq1$,}\\
						       & \text{  multiple roots satisfy $|\zeta(\mu)|<1$} \end{split}
\right\}.
\end{equation}
Solving for $\mu=\Delta t \lambda$ in (\ref{eq:SBDF5}), we get that our method is stable for 
\begin{equation}
\mu=\frac{\frac{-3}{2}\zeta^2 +2\zeta -\frac{1}{2}}{1-2\zeta}, \qquad\mbox{with}\quad|\zeta|\leq 1.
\end{equation}
We take $\zeta = e^{i\theta}$ for $0\leq \theta \leq 2\pi$ to find the contour of the stability region:
\begin{equation}
\mu = \frac{\frac{-3}{2}e^{2i\theta} +2e^{i\theta} -\frac{1}{2}}{1-2e^{i\theta}},\quad\theta\in [0,2\pi].
\end{equation}
The contour of the stability region is shown in Figure  \ref{fig:stabsbdf}. As with the previous methods, the most restrictive eigenvalue will be the most negative. This occurs for a real eigenvalue, hence $-4/3\leq \lambda \Delta t \leq 0$ for stability. The critical time-step is then $\Delta t_{\mathrm{theo}}^*=\dfrac{-4}{3\lambda_{\mathrm{min}}}$ and its value is shown in Table \ref{tab:crittheo}  for the BR model, in Table \ref{tab:crittheoMS} for the MS model and in Table \ref{tab:crittheoTNNP} for the TNNP model.

\begin{remark}[Justification for the choice of $\omega$]

We rewrite equation (\ref{eq:SBDF2}) as a one-step recursive equation:
\begin{equation}
\begin{bmatrix} Y^{n+1} \\ Y^n \end{bmatrix} = 
\begin{bmatrix}
(\frac{3}{2}I_d-A_\omega)^{-1}(2I_d+2\Delta t J) & (\frac{3}{2}I_d-A_\omega)^{-1}(-\frac{1}{2}I_d-\Delta t J)\\
I_d & 0
\end{bmatrix}
\begin{bmatrix} Y^n \\ Y^{n-1} \end{bmatrix}.
\end{equation}
Let 
\begin{equation}
\displaystyle R(\Delta t, h, \omega)=\begin{bmatrix}
(\frac{3}{2}I_d-A_\omega)^{-1}(2I_d+2\Delta t J) & (\frac{3}{2}I_d-A_\omega)^{-1}(-\frac{1}{2}I_d-\Delta t J)\\
I_d & 0
\end{bmatrix},
\end{equation}
be the stability function for the SBDF2 method. We can factor out from $R$ the matrix
\begin{equation}\begin{bmatrix}
(I_d-\frac{2}{3}A_\omega)^{-1} & 0\\
0 & I_d
\end{bmatrix},
\end{equation}
which has a spectral radius of 1. This leads to the stability condition $\rho(\tilde{R})\le 1$, where
\begin{equation}
\displaystyle \tilde{R}(\Delta t, h, \omega)=\begin{bmatrix}
\frac{2}{3}(2I_d+2\Delta t J) & \frac{2}{3}(-\frac{1}{2}I_d-\Delta t J)\\
I_d & 0
\end{bmatrix}.
\end{equation}
The matrix $\tilde{R}$ appears in the one-step formulation of the difference equation (\ref{eq:SBDF3}).
\end{remark}
\subsubsection{Rush-Larsen} 
We now look at the stability of the RL-FBE and the RL-CNAB schemes. These schemes differ from the methods previously studied in that the differential equations for the gating variables are solved using a different method than for the transmembrane potential and the concentrations. The major advantage of using Rush-Larsen methods is to be able to use large time-steps compared to more classical methods. In fact, it can be observed that in the case of most problems in electrophysiology, the region of stability of RL methods of the form (\ref{eq:RL2}) covers the entire negative half plane \cite{PerVen2009}, i.e. A-stability. The stability of the scheme then depends on the methods used to solve the differential equations for the transmembrane potential and the concentrations. The following discussion on the stability of the Rush-Larsen methods uses heuristic arguments that provide critical time-steps relatively close to those in numerical tests, but this derivation cannot yet be formalized.

We use a different scheme for the differential equations of the gating variables and the concentrations. Therefore, we look at the problem in the form (\ref{eq:nondim01})-(\ref{eq:nondim03}). As done previously, we linearize around some constant state $(\tilde{u},\tilde{\mathbf{v}}, \tilde{\mathbf{X}})=\tilde{Y}$ and take the Fourier transform. We apply the scheme (\ref{eq:RLscheme}) to the linearized problem. The scheme is identical to the Forward-Backward Euler method for the transmembrane potential and concentrations. We get
\begin{equation}\label{eq:RLU}
\left( 1- \sigma\dfrac{2 \cos \omega h -2}{h^2}\right) U_\omega^{n+1}=U_\omega^{n}+\Delta t 
\begin{bmatrix}  -\dfrac{\partial I}{\partial u}(\tilde{Y})  & -\dfrac{\partial I}{\partial \mathbf{v}}(\tilde{Y}) & -\dfrac{\partial I}{\partial \mathbf{X}}(\tilde{Y}) \end{bmatrix} 
\begin{bmatrix}	U_\omega^{n}\\	V_\omega^{n}\\	X_\omega^{n}\end{bmatrix},
\end{equation}
\begin{equation}\label{eq:RLX}
X_\omega^{n+1}=X_\omega^{n}+\Delta t
 \begin{bmatrix}  \dfrac{\partial \mathbf{g}}{\partial u}(\tilde{Y})  & \dfrac{\partial \mathbf{g}}{\partial \mathbf{v}}(\tilde{Y}) & \dfrac{\partial \mathbf{g}}{\partial \mathbf{X}}(\tilde{Y}) \end{bmatrix} 
\begin{bmatrix}	U_\omega^{n}\\	V_\omega^{n}\\	X_\omega^{n}\end{bmatrix}.
\end{equation}
We solve the ODEs for the gating variables with the RL method, which in this case is an exponential integrator method. As opposed to previously studied methods, it is not a linear multistep method. Because the Rush-Larsen method is considered A-stable when applied to the gating equations \cite{PerVen2009}, we assume for the sake of analysis that we have an A-stable one-step linear method for $V$, which can be written as
\begin{equation}\label{eq:RLonestep}
V_\omega^{n+1}=R(\Delta t,\omega)V_\omega^{n},
\end{equation}
where $\rho(R(\Delta t,\omega))\leq 1$, for any $\Delta t>0$. We then consider the most restrictive case where $\rho(R(\Delta t,\omega)) = 1$. This occurs when $R(\Delta t,\omega)) = I$, which is equivalent to setting
\begin{equation}\label{eq:RLV}
V_\omega^{n+1}=V_\omega^{n}.
\end{equation}
As for the previous semi-implicit methods, we take $\omega = 0$ because it is the choice that is most restrictive for stability. Combining (\ref{eq:RLU})-(\ref{eq:RLV}), we have
\begin{equation}\label{eq:RLstab}
Y^{n+1}=(I+\Delta t J_{RL} ) Y^n,
\end{equation}
where
\begin{equation}\label{eq:RLjac}
J_{RL}=
\begin{bmatrix}
-\dfrac{\partial I}{\partial u}(\tilde{Y})  & -\dfrac{\partial I}{\partial \mathbf{v}}(\tilde{Y}) & -\dfrac{\partial I}{\partial \mathbf{X}}(\tilde{Y})\\[0.5em] 
0 & 0& 0\\
\dfrac{\partial \mathbf{g}}{\partial u}(\tilde{Y})  & \dfrac{\partial \mathbf{g}}{\partial \mathbf{v}}(\tilde{Y}) & \dfrac{\partial \mathbf{g}}{\partial \mathbf{X}}(\tilde{Y})
\end{bmatrix}.
\end{equation}
As we did before for the Euler methods, the stability condition for this scheme is $\Delta t \leq -2/\lambda_{\mathrm{min}}$, where $\lambda_{\mathrm{min}}$ is the most negative eigenvalue of $J_{RL}$. The critical time-step, $\Delta t^*_{\mathrm{theo}}=-2/\lambda_{\mathrm{min}}$ is shown in Table \ref{tab:crittheo} for the BR model and in Table \ref{tab:crittheoTNNP} for the TNNP model.

Similarly, for the RL-CNAB method, we get from the stability analysis of the Adams Bashforth method, $\Delta t^*_{\mathrm{theo}}=-1/\lambda_{\mathrm{min}}$.

\subsubsection{Deferred Correction} 

Due to the more complex nature of the third-order deferred correction scheme, we cannot easily find a stability condition using absolute stability analysis. However, the numerically observed critical time-step is very close to the one from the Forward-Backward Euler method.

\subsection{Numerical Results}\label{sec:numres}
In this section, we compare the critical time-steps obtained through our absolute stability analysis from the previous section to those observed numerically in the case of the BR, MS and TNNP models. All parameters used and the details on the ionic models are given in \cite{Roy2015}. 

For the 1D case, we use a spatial domain of length $\SI{100}{\centi\metre}$, discretized by equally spaced nodes $x_i = ih$, $h=100/N$ or its equivalent for the nondimensionalized MS model. We used a final time $T$ of $\SI{400}{\milli\second}$ for the BR model, $\SI{350}{\milli\second}$ for the MS model and $\SI{300}{\milli\second}$ for the TNNP model. 
The applied stimulation current is the $C^\infty$ function given by
\begin{equation}
I_{app}(\mathbf{x},t) = \begin{cases} 	50 \exp \left(1-\dfrac{1}{1-(t-1.5)^2}\right) \exp \left(1-\dfrac{1}{1-4(x-0.5)^2}\right), &\begin{split} \text{if } 0<&x<1, \\ 0.5<&t<2.5, \end{split}\\
					0, & \text{otherwise}.
		   \end{cases}
\end{equation}
All computations for the 1D case were made with MATLAB.

We also tested the stability of the BR model in the 2D case using the code CHORAL \cite{Choral}. For the 2D simulations, we used a $1\si{\centi\metre}\times 1\si{\centi\metre}$ square domain discretized with 3432 grid points and a final time $T$ of  $\SI{16}{\milli\second}$. The applied stimulation current in the 2D case is a $C^1$ function of $x$ and $t$ given by
\begin{equation}\label{Istim2D}
I_{app}(\mathbf{x},t) = \begin{cases} 	50 \dfrac{1 + \cos (\pi r)}{2}\dfrac{1 + \cos (\pi \tau)}{2},  &\text{if } 0\leq r \leq 1, \; 0 \leq \tau \leq 1, \\
					0, & \text{otherwise},
		   \end{cases}
\end{equation}
where $r = |\mathbf{x} - \mathbf{x_0}|/r_0$ with $\mathbf{x_0}$ the point at the center of a circular simulation zone and $r_0$, the radius of this zone. Similarly, $\tau = (t-t_0)/\tau_0$ sets a stimulation interval in time starting at $t=t_0$ with duration $\tau_0$. We set $t_0=5\si{\milli\second}$,  $\tau_0=1.5\si{\milli\second}$ and $r_0=0.125\si{\centi\metre}$.

For each method, the theoretical critical time-step $\Delta t^*_{\mathrm{theo}}$ is determined by the stability conditions from the previous section and the most negative eigenvalue of the Jacobian of the ionic model used, $\lambda_{\mathrm{min}}$, which is calculated on the domain as explained in the previous section. The numerically observed critical time-step $\Delta t^*$ is the largest possible time-step for which the numerical solution remains bounded. When using a relatively large time-step with the more stable Rush-Larsen methods, the solution remains bounded but its shape degenerates and does not represent the real shape of the wave. Therefore, for the RL methods the value of $\Delta t^*$ is taken as the largest possible time-step for which the potential wave does not degenerate. 

\begin{remark}[Computing the Jacobian for discontinuous models]

The functions $I$, $F$ and $G$ of the BR model are all continuously differentiable with respect to each variable. We thus found the Jacobian analytically. Other ionic models have discontinuities which prevent us from defining a Jacobian for discontinuity points. To avoid deriving expressions on both sides of the discontinuities, we decided to approximate the Jacobian numerically, for e.g. by using MATLAB's \texttt{numjac} function. Several functions in the ODEs for the gating variables in the TNNP model have discontinuities at $u=-40\si{\milli\volt}$ and the function in the ODE of the gating variable of the MS model has a discontinuity at $u=u_{gate}$. These are null sets of the phase space and in general we will not observe these discontinuities when approximating the Jacobian or solving the differential equations numerically. However, if we approximate the Jacobian at the discontinuity point or close enough to be below the tolerance for \texttt{numjac}, we obtain extremely large eigenvalues. These discontinuities might become a problem for extremely small time-steps because there could be values of $u$ very close to the singularities. This has not been a problem for our simulations. 
\end{remark}
\subsubsection{Beeler-Reuter}

\begin{table}[!htb]
\caption{Size of $\Delta t^*$ for the numerical methods used with BR model}\label{tab:critnum}
\sisetup{round-mode=figures, round-precision=5,table-parse-only}
\centering
\begin{tabular}{| l | | S | S | S |}
\hline
Methods & $h=0.0625$ & $h=0.03125$ & $h=0.015625$\\
\hline \hline
2\textsuperscript{nd} Order SBDF 	&0.0168477803 	&0.0168477803 	& 0.0168477803\\
\hline
Strang Splitting (CN-RK4)		&0.0714796283 	&0.0714796283 	&0.0714796283 \\
\hline
Strang Splitting (CN-RK2)		&0.0502891627	&0.0502891627	&0.0502891627 \\
\hline
Forward Euler				&0.025372661	&0.0200944439 	& 0.0050639321\\
\hline
Forward-Backward Euler 		&  0.025372661 	& 0.025372661 	& $0.025372661$\\
\hline
RL-CNAB 					&0.2352941176	& 0.2352941176 	& 0.2352941176 \\
\hline
RL-FBE				&$>0.8$		& $>0.8$ 		& $>0.8$ \\
\hline
DC3				&0.0244648318	&0.0244648318		&0.0244648318 \\
\hline
\end{tabular}

 \vspace*{\floatsep}

\caption{Size of $\Delta t^*_{\mathrm{theo}}$ for the numerical methods used with BR model}\label{tab:crittheo}
\sisetup{round-mode=figures, round-precision=5,table-parse-only}
\centering
\begin{tabular}{| l | | S | S | S |}
\hline
Methods & $h=0.0625$ & $h=0.03125$ & $h=0.015625$\\
\hline \hline
2\textsuperscript{nd} Order SBDF 	&0.0163035284 	& 0.0163035284	&0.0163035284\\
\hline
Strang Splitting (CN-RK4)		&0.068115169	&0.068115169 	&0.068115169 \\
\hline
Strang Splitting (CN-RK2)		&0.0489105852	&0.0489105852	&0.0489105852 \\
\hline
Forward Euler 				&0.0244552578	&0.0202383975 	&0.0050661773\\
\hline
Forward-Backward Euler 		&  0.0244552926	& 0.0244552926	&0.0244552926\\
\hline
RL-CNAB 					&0.4234156457 	& 0.4234156457 	& 0.4234156457\\
\hline
RL-FBE			&0.8468312914	&0.8468312914 	& 0.8468312914 \\
\hline
\end{tabular}

 \vspace*{\floatsep}

\caption{Size of $\Delta t^*$ for the numerical methods used with BR model in 2D}\label{tab:critnum2D}
\sisetup{round-mode=figures, round-precision=5,table-parse-only}
\centering
\begin{tabular}{| c| | S |}
\hline
Methods & 	$\Delta t^*$ \\
\hline \hline
2\textsuperscript{nd} Order SBDF 	&0.016131	\\
\hline
Strang Splitting (CN-RK4)		&	0.059566 \\
\hline
Strang Splitting (CN-RK2)		&	0.044933 \\
\hline
Forward-Backward Euler 		&  0.024242	\\
\hline
RL-CNAB 					& 	0.2 \\
\hline
\end{tabular}
\end{table}

The most negative eigenvalue of the Jacobian for the Beeler-Reuter model is $\lambda_{\mathrm{min}}=-81.782$.
Note that one can find the most negative eigenvalues of the Jacobian of 37 different ionic models in ~\cite{SpitDean2010,Marsh2013}, where the authors obtained a value of -82.0 for the BR model. This small difference is most likely due to a distinct applied stimulation current or initial conditions. The theoretical critical time-steps for methods studied above are shown in Table \ref{tab:crittheo} and the numerically observed critical time-steps in Table \ref{tab:critnum} for the 1D case and in Table \ref{tab:critnum2D} for the 2D case. As expected, the critical time-steps for the semi-implicit methods are independent of $h$. 

We see that for all methods except the RL methods, the critical time-steps obtained numerically are very close to those obtained through the absolute stability analysis. In the case of the RL methods, the critical time-steps are similar for RL-FBE, but $\Delta t^*$ is approximately half of $\Delta t^*_{\mathrm{theo}}$ for RL-CNAB. This is likely a consequence of the difficulty to identify the critical time step in numerical tests.

For the Strang splitting methods, we observe smaller $\Delta t^*$ in 2D compared to their 1D equivalent. 
Otherwise, the analysis done in the last section for the 1D monodomain model seems to apply for the semi-implicit methods in the 2D case. Due to the more complex nature of the Finite Element method for higher dimensions, the analysis of explicit methods such as FE depends on the nature of the mesh used (uniform vs non-uniform). 

We observe that the critical time-steps of the FE and FBE methods are initially the same for $h=0.0625$, but as $h$ gets smaller FE becomes less stable with a critical time-step of order $h^2$. This indicates that the use of explicit methods to solve the BR model is only justifiable for very coarse meshes.

The RL methods are the most stable of all the methods studied. The RL-CNAB has a $\Delta t^*$ more than three times larger than the one of the next most stable method, CN-RK4. This method has a $\Delta t^*$ slightly larger than the one for CN-RK2, which reflects the fact that the stability region of CN-RK2 is included in the stability region of CN-RK4. The $\Delta t^*$ for the SBDF2 method is three to four times smaller than the ones for the Strang splitting methods and $50\%$ smaller than the one for FBE. These relations between the $\Delta t^*$ of the linear multistep methods, i.e. excluding the RL methods, are the same for all the ionic models used. 

As mentioned in the last section, the critical time-step for the DC3 method is very close to the value for FBE: it is slightly smaller.

\subsubsection{Mitchell-Schaeffer}

\begin{table}[!htb]
\caption{Size of $\Delta t^*$ for the numerical methods used with MS model}\label{tab:critnumMS}
\sisetup{round-mode=figures, round-precision=5,table-parse-only}
\centering
\begin{tabular}{| l | | S | S | S |}
\hline
Methods & $h=1$ & $h=0.5$ & $h=0.25$\\
\hline \hline
2\textsuperscript{nd} Order SBDF 	&0.5263157895 	&0.5263157895 	& 0.5263157895\\
\hline
Strang Splitting (CN-RK4)		&1.8518518519	&1.8518518519   &1.8518518519 \\
\hline
Strang Splitting (CN-RK2)		&1.5217391304	&1.5217391304	&1.5217391304	 \\
\hline
Forward Euler				&0.1226993865	&0.0358744395	& 0.0089405454\\
\hline
Forward-Backward Euler 		& 0.7692307692 	& 0.7692307692	& 0.7692307692\\
\hline
DC3 		& 0.7692307692 	& 0.7692307692	& 0.7692307692\\
\hline
\end{tabular}
 \vspace*{\floatsep}

\caption{Size of $\Delta t^*_{\mathrm{theo}}$ for the numerical methods used with MS model}\label{tab:crittheoMS}
\sisetup{round-mode=figures, round-precision=5,table-parse-only}
\centering
\begin{tabular}{| l | | S | S | S |}
\hline
Methods & $h=1$ & $h=0.5$ & $h=0.25$\\
\hline \hline
2\textsuperscript{nd} Order SBDF 	&0.5002868518 	&0.5002868518	&0.5002868518\\
\hline
Strang Splitting (CN-RK4)		&2.0901686224	&2.0901686224	&2.0901686224 \\
\hline
Strang Splitting (CN-RK2)		&1.5008605555	&1.5008605555	&1.5008605555 \\
\hline
Forward Euler 				&0.111790333	&0.033672183	&0.0088728184\\
\hline
Forward-Backward Euler 		&  0.7504302777	& 0.7504302777	&0.7504302777\\
\hline
\end{tabular}
\end{table}

For the MS model, the most negative eigenvalue obtained is $\lambda_{\mathrm{min}} =-2.6651$. The theoretical critical time-steps for the methods studied are shown in Table \ref{tab:crittheoMS} and the numerically observed critical time-steps in Table \ref{tab:critnumMS}. We see that for all methods, $\Delta t^*$ is very close to $\Delta t^*_{\mathrm{theo}}$. As expected, the critical time-steps for the semi-implicit implicit methods are independent of the space step $h$. 
Because the MS model is not very stiff, the stability of the FE method depends on the size of $h$ even for relatively large space step $h$. This indicates the necessity for taking the diffusion implicitly when solving less stiff models such as the MS model.

We also observe that the CN-RK4 method requires a slightly smaller time-step than $\Delta t^*_{\mathrm{theo}}$. This is most likely resulting from the oscillations caused by the use of the Crank-Nicolson method \cite{torabi2014}. These oscillations can be observed for large time-steps for both Strang splitting methods.

For the MS model, the Strang splitting methods are the most stable of all the methods studied. They have a $\Delta t^*$ three to four times larger than for the SBDF2 method. 

The DC3 and FBE methods have the same numerically observed critical time-step.

\subsubsection{ten Tuscher-Noble-Noble-Panfilov}
\begin{table}[!htb]
\caption{Size of $\Delta t^*$ for the numerical methods used with TNNP model}\label{tab:critnumTNNP}
\sisetup{round-mode=figures, round-precision=5,table-parse-only}
\centering
\begin{tabular}{| l | | S | S | S |}
\hline
Methods & $h=0.0625$ & $h=0.03125$ & $h=0.015625$\\
\hline \hline
2\textsuperscript{nd} Order SBDF 	&0.0011348205 	&0.0011348205 	& 0.0011348205\\
\hline
Strang Splitting (CN-RK4)		&0.0048362136 	&0.0048362136 	&0.0048362136 \\
\hline
Strang Splitting (CN-RK2)		&0.0034562212	&0.0034562212	&0.0034562212 \\
\hline
Forward Euler				&0.0017044196	&0.0017044196 	& 0.0017044196\\
\hline
Forward-Backward Euler 		&  0.0017044196 	& 0.0017044196 	& $0.0017044196$\\
\hline
RL-CNAB					&0.0913798355	& 0.0913798355 	& 0.0913798355 \\
\hline
RL-FBE			&$>0.6$		& $>0.6$ 		& $>0.6$ \\
\hline
DC3				& 0.0016869096		& 0.0016869096	& 0.0016869096 \\
\hline
\end{tabular}

 \vspace*{\floatsep}
 
\caption{Size of $\Delta t^*_{\mathrm{theo}}$ for the numerical methods used with TNNP model}\label{tab:crittheoTNNP}
\sisetup{round-mode=figures, round-precision=5,table-parse-only}
\centering
\begin{tabular}{| l | | S | S | S |}
\hline
Methods 					& $h=0.0625$ 		& $h=0.03125$ 		& $h=0.015625$\\
\hline \hline
2\textsuperscript{nd} Order SBDF 	&0.0011188798		& 0.0011188798		&0.0011188798\\
\hline
Strang Splitting (CN-RK4)		&0.0046746132464	&0.0046746132464	&0.0046746132464\\
\hline
Strang Splitting (CN-RK2)		&0.003356639535	&0.003356639535	&0.003356639535 \\
\hline
Forward Euler				&  0.0016783198		&   0.0016783198		&  0.0016783198\\
\hline
Forward-Backward Euler 		&   0.0016783198		&   0.0016783198		&  0.0016783198\\
\hline
RL-CNAB 					& 0.2061158342		& 0.2061158342 		& 0.2061158342\\
\hline
RL-FBE 				&0.4122316684		& 0.4122316684 		& 0.4122316684 \\
\hline
\end{tabular}

 \vspace*{\floatsep}
 
\caption{Size of $\Delta t^*$ for the Forward Euler's method with TNNP model}\label{tab:FEtnnp}
\sisetup{round-mode=figures, round-precision=5,table-parse-only}
\centering
\begin{tabular}{| c| | S | S |}
\hline
$h$ 		&$\Delta t^*$  	& $\Delta t_{\mathrm{theo}}^*$ \\	\hline\hline
0.015625 	&0.0017043518	&	0.0016783198	\\\hline
0.0078125	&0.0012677485	&	0.0012668534	\\\hline
0.00390625 &0.0003168233	&	0.0003167574	\\\hline

\end{tabular}
\end{table}

\begin{sloppypar}For the TNNP model, the most negative eigenvalue obtained is $\lambda_{\mathrm{min}} = -1191.7$. The value given in ~\cite{SpitDean2010,Marsh2013} is -1170. Again, this difference is most likely due to a distinct applied stimulation current or initial conditions. The theoretical critical time-steps for the methods studied are shown in Table \ref{tab:crittheoTNNP} and the numerically observed critical time-steps in Table \ref{tab:critnumTNNP}. 

As for the BR model, we see that for all methods except the RL methods, the critical time-steps obtained numerically are very close to those obtained through absolute stability analysis. As expected, the critical time-steps for the semi-implicit methods are independent of the space step $h$. 

In the case of the RL methods, $\Delta t^*$ is approximately half of $\Delta t^*_{\mathrm{theo}}$ for RL-CNAB and two thirds for RL-FBE.\end{sloppypar}
For the Forward Euler scheme, $\Delta t^*$ is the same for $h=0.0625, 0.03125,\allowbreak 0.015625$. For smaller values of $h$, we begin to see a dependence of $\Delta t^*$ on $h^2$, as shown in Table \ref{tab:FEtnnp}. These results indicate that for very stiff models such as the TNNP model, the use of fully explicit methods could still be acceptable, except on finer spatial meshes. 

The RL methods are the most stable of all the methods studied. The RL-CNAB has a $\Delta t^*$ almost twenty times larger than the next most stable method, CN-RK4.

As for the BR model, the DC3 method has a $\Delta t^*$ slightly smaller than the one for FBE.

As expected from the analysis of the last sections, for all ionic models, only the fully explicit Forward Euler's method has a critical time-step depending on the size of $h$. Indeed, as $h$ gets smaller we see that $\Delta t^*$ is eventually proportional to $h^2$. Looking at the Jacobian for large $h$, we see that the minimum eigenvalue comes from the gating variables, in the case of the models studied. For the FE method, small values of $h$ are required for the eigenvalue coming from the diffusion term in the monodomain equation to become the most negative. This leads to the conclusion that for very stiff models such as the TNNP model, taking the diffusion implicitly is only necessary for very fine meshes.

For all models, the semi-implicit method which requires the smallest time-step is the SBDF2 method.

\section{Accuracy of the numerical methods}\label{chap:convergence}

In this section, we investigate the accuracy of the different time-stepping methods when solving the monodomain model coupled with the three different ionic models studied. We start by conducting a convergence test for each method to verify if they have the correct rate of convergence. The accuracy of the methods will also be compared relatively to the size of the time-step used. Afterwards, we will compare the accuracy of the methods relatively to the computational time needed to run the simulations.

\subsection{Convergence tests}

We now study the convergence of the different time-stepping methods for a given spatial mesh. We split the error in space and time given that their order may not be the same and study the error in time only. We have
\begin{equation}
\|u - u_{h,\Delta t} \| \leq \| u - u_h \| + \|u_h -  u_{h,\Delta t} \| = \mathcal{O} (h^p + \Delta t^q),
\end{equation}
where $u$ is the exact solution for the transmembrane potential, $u_h$ is the solution of the semi-discretized problem in space and $u_{h,\Delta t}$ is the solution of the fully discretized problem. Since we want to study the convergence of the time-stepping methods, we will only consider the error between the solutions of the semi-discretized and fully discretized problems. A consequence of considering only the error in time is that relative error levels requested in our test cases will be smaller (sometimes much smaller) than relative error levels in space and time usually found in the cardiac electrophysiology literature.

We test convergence with respect to the following errors:
\begin{align}
\label{eq:L2error}
\begin{split}
e_{L^2}  =  \| u_{h,\Delta t}(T)-u_h(T) \|_{L^2}, &\qquad e_{H^1} = | u_{h,\Delta t}(T)-u_h(T) |_{H^1},
\\
e_c = |c_{h,\Delta t}-c_h|,&\qquad e_{T_1} = |T_{1,h,\Delta t}-T_{1,h}|,
\end{split}
\end{align}
where $u_h$ is a reference semi-discretized solution for the transmembrane potential, which is calculated using a very small time-step. The norms used are discrete approximations of the $L^2(\Omega)$ norm and $H^1(\Omega)$ seminorm using Simpson's rule. The solutions $u_{h,\Delta t}$ are calculated with the same spatial mesh and at the same final time $T$ as the reference solution, but using larger values for $\Delta t$.  We denote by $c_{h,\Delta t}$ the wave velocity, and by $T_{1,h,\Delta t}$ the depolarization time, i.e. the time at which a given point of the domain reaches a given super-threshold value of the transmembrane potential. The wave velocity and depolarization time of the reference solution $u_h$ are denoted by $c_h$ and $T_{1,h}$, respectively. The wave velocity is defined by $c= (x_2-x_1)/(T_2-T_1)$, where $T_i$ is the time when the depolarization front of the potential wave passes through chosen nodes $x_i$, $i=1,2$. The wave front passes through a point $x_i$ during the time-step from $t_n$ to $t_{n+1} = t_n + \Delta t$ if $u(x_i,t_n)<\hat{u}$, but $u(x_i,t_{n+1})\geq\hat{u}$ for some chosen value $\hat{u}$ on the wave front.
For better approximations of $T_i$, we use linear interpolation and define
\begin{equation}
T_i = t_n + \Delta t \frac{\hat{u} - u(x_i,t_n)}{u(x_i,t_{n+1}) - u(x_i,t_n)}.
\end{equation}
Assuming that the error is proportional to $\Delta t^\alpha$, the estimated convergence rate is calculated with
\begin{equation}\label{eq:rate}
\alpha = \frac{\log (|e_1/e_2|)}{\log (\Delta t_1/\Delta t_2)},
\end{equation}
where $\Delta t_1$ and $\Delta t_2$ are consecutive time-steps in a sequence of decreasing time-steps, and $e_1$ and $e_2$ are the corresponding errors. All solutions are calculated using the same stimulation current and parameters from Section \ref{sec:numres}. The largest $\Delta t$ used for each method is taken close to the critical time-step found in the previous section, except for the RL methods where we use a starting time-step similar to the other methods, that is when the RL methods start to show their asymptotic behaviour. 

The order of convergence in the $L^2$ norm and $H^1$ seminorm are expected to be the same because the numerical solutions $u_{h,\Delta t}$ are in the same finite-dimensional space and thus the equivalence of norms applies.

\subsubsection{Beeler-Reuter with 1D Monodomain}
\begin{figure}[htpb!]
\centering
\begin{subfigure}[b]{0.31\textwidth}
\includegraphics[width=\linewidth]{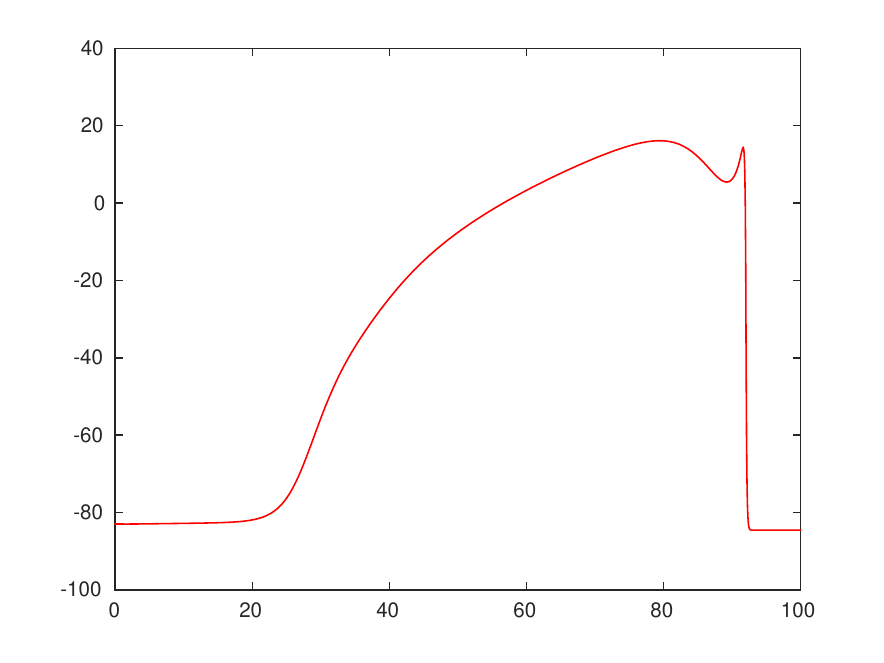}

  \caption{BR model: $u$ at time $T=\SI{400}{\milli\second}$, plotted for $x\in [0,100]$.} \label{fig:refsol}
  \end{subfigure}
   \hfill
  \begin{subfigure}[b]{0.31\textwidth}
\includegraphics[width=\linewidth]{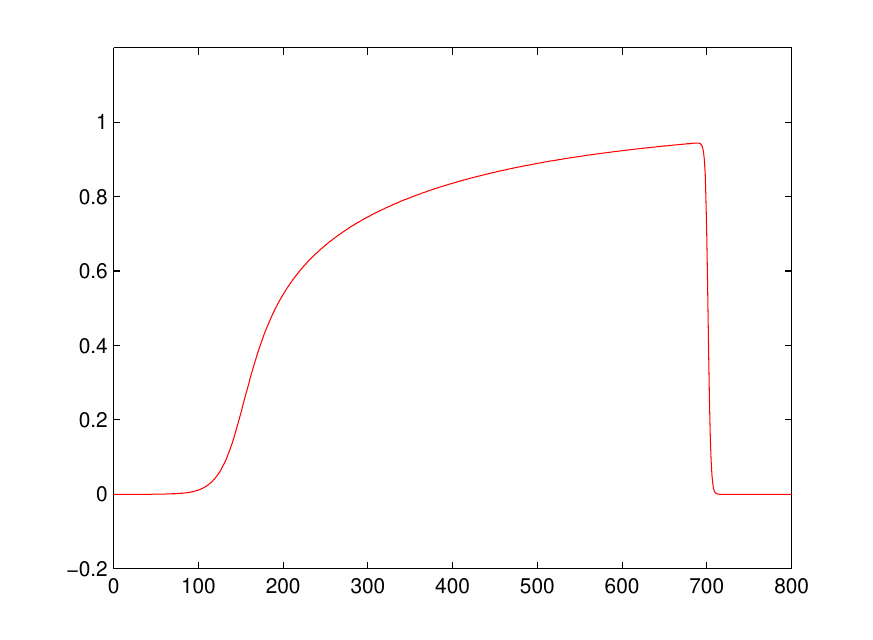}
 
  \caption{MS model: $u$ at time $T=\SI{350}{\milli\second}$, plotted for $x\in [0,800]$.} \label{fig:refsolMS}
\end{subfigure}
 \hfill
\begin{subfigure}[b]{0.31\textwidth}
\includegraphics[width=\linewidth]{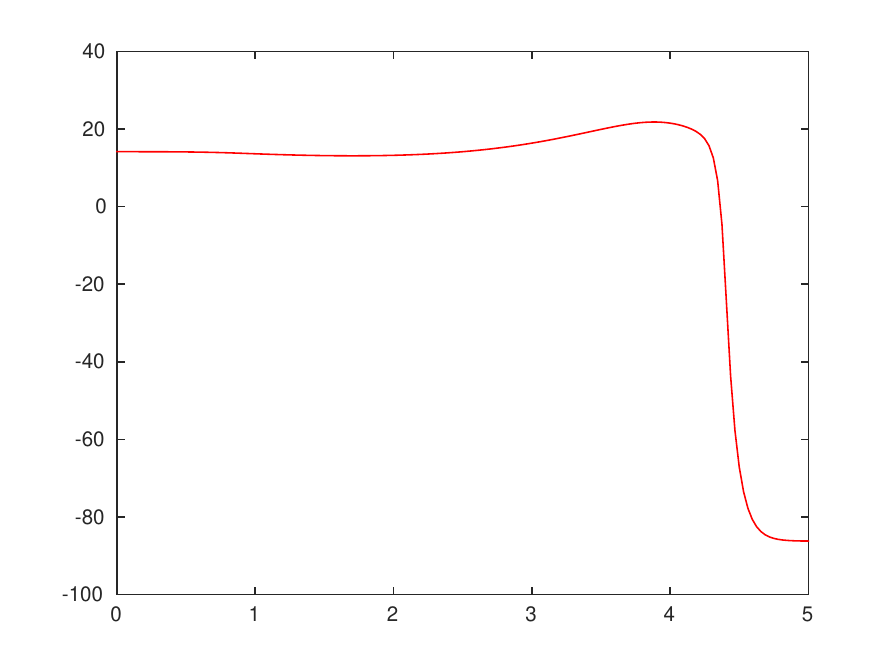}
 
  \caption{TNNP model: $u$ at time $T=\SI{12}{\milli\second}$, plotted for $x\in [0,5]$.} \label{fig:refsoltnnp}
\end{subfigure}

\caption{Reference solutions for the transmembrane potential}
  
\end{figure}

The reference solution $u_h$ for the BR model is obtained using SBDF2 on a domain of length $100 \;cm$ discretized in space with $1600$ nodes and computed at time $T = \SI{400}{\milli \second}$ with $\Delta t = \SI{1/125000}{\milli \second}$. The graph of the reference solution is shown in Figure \ref{fig:refsol}. 

\begin{figure}[!htb]
\centering
\begin{subfigure}[b]{0.49\textwidth}
\includegraphics[width=\linewidth]{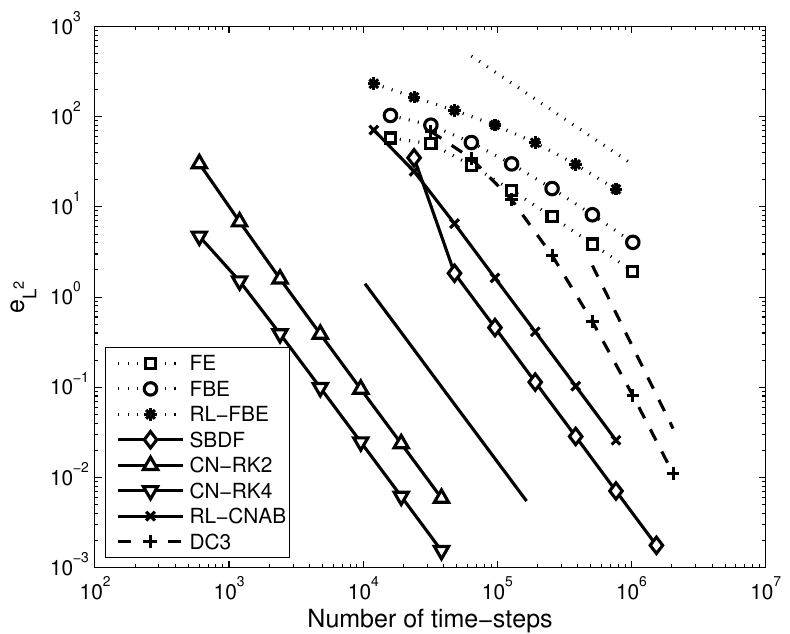}

  \caption{Error in $L^2(\Omega)$ norm} \label{fig:BRL2}
  \end{subfigure}
   \hfill
  \begin{subfigure}[b]{0.49\textwidth}
\includegraphics[width=\linewidth]{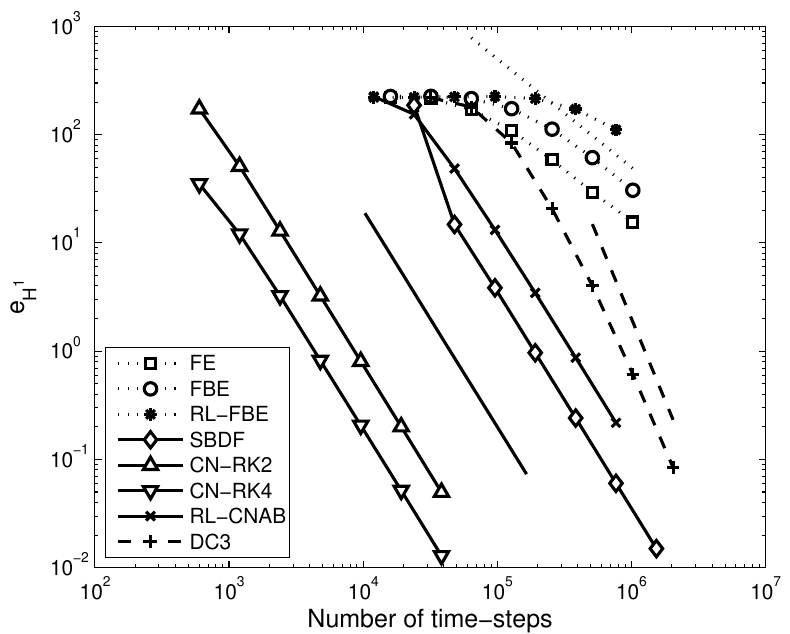}
 
  \caption{Error in $H^1(\Omega)$ seminorm} \label{fig:BRH1}
\end{subfigure}
 \vskip\baselineskip
 \begin{subfigure}[b]{0.49\textwidth}
\includegraphics[width=\linewidth]{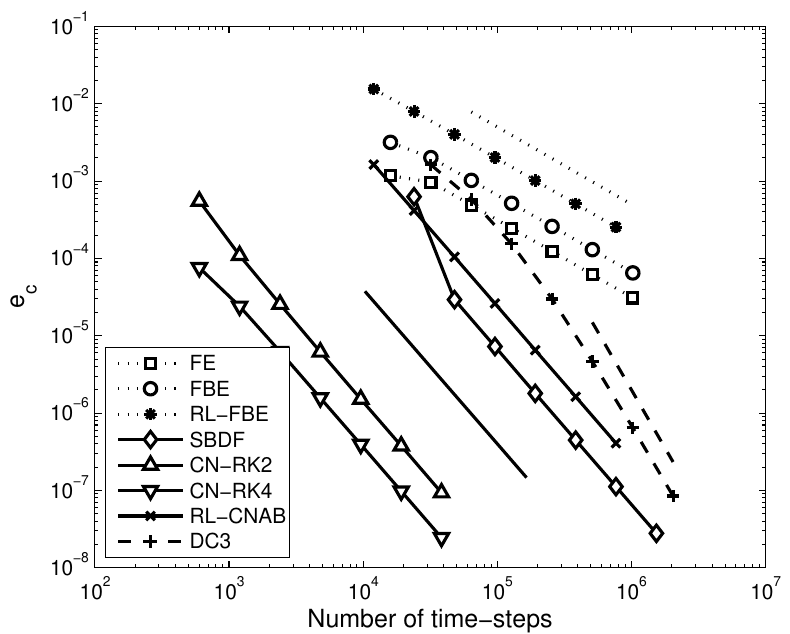}

  \caption{Error of the wave velocity $c$} \label{fig:BRc}
  \end{subfigure}
   \hfill
  \begin{subfigure}[b]{0.49\textwidth}
\includegraphics[width=\linewidth]{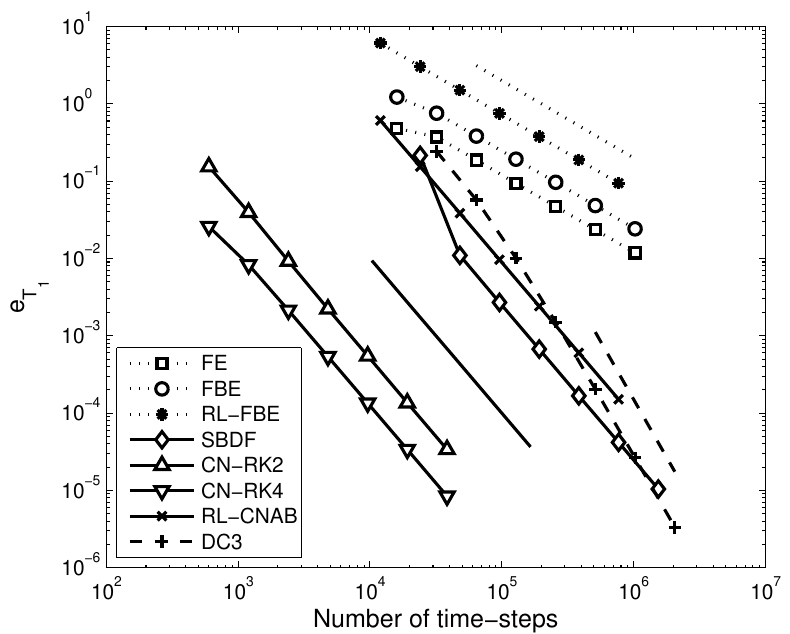}
 
  \caption{Error of the depolarization time $T_1$} \label{fig:BRT}
\end{subfigure}

\caption{Errors for the BR model at fixed $h=1/16$. The order of convergence can be observed by comparing with lines of negative slope one (dotted), two (solid), and three (dashed).}\label{fig:BRconv}
  
\end{figure}

To determine $c$ and $T_1$, we take $\hat{u} = \SI{-30}{\milli \volt}$, $x_1 = \SI{20}{\centi\meter}$ and $x_2 = \SI{50}{\centi\meter}$. The results of the convergence tests for the 1D BR model are shown in Figure \ref{fig:BRconv}. We illustrate how the  errors defined in (\ref{eq:L2error}) vary as the size of the time-step is decreased. We observe that $T_1$ and $c$ reach their asymptotic rate of convergence faster than the $L^2$ norm and $H^1$ seminorm of the error.

For the 1D BR model, all of the methods studied showed their expected asymptotic order of convergence. For the first-order methods, we observe that for the same $\Delta t$, FE is almost twice as accurate as its semi-implicit version, FBE, and both methods are significantly more accurate than RL-FBE. It takes very small values of $\Delta t$ for RL-FBE to reach its asymptotic order of convergence.

For the second-order methods, we observe that for the same $\Delta t$, CN-RK4 is about ten times more accurate than CN-RK2, which in turn is about two times more accurate than SBDF2. This last method is two times more accurate than RL-CNAB. For the DC3 method, it takes very small values of $\Delta t$ before it exhibits its correct order of convergence. Even for $\Delta t = 1/5120$, the error is still larger than what we would obtain by using a second-order method. We noticed that this unexpected reduced accuracy of DC3 sometimes results when the stimulation current is applied close to the boundary of the domain.

\subsubsection{Mitchell-Schaeffer with 1D Monodomain}
The reference solution for the MS model is computed using SBDF2 on a domain of length $800$ discretized in space with $800$ nodes and computed at time $350$ with $\Delta t = 7/60000$. The plot of the reference solution is shown in Figure \ref{fig:refsolMS}.

\begin{figure}[!htb]
\centering
\begin{subfigure}[b]{0.49\textwidth}
\includegraphics[width=\linewidth]{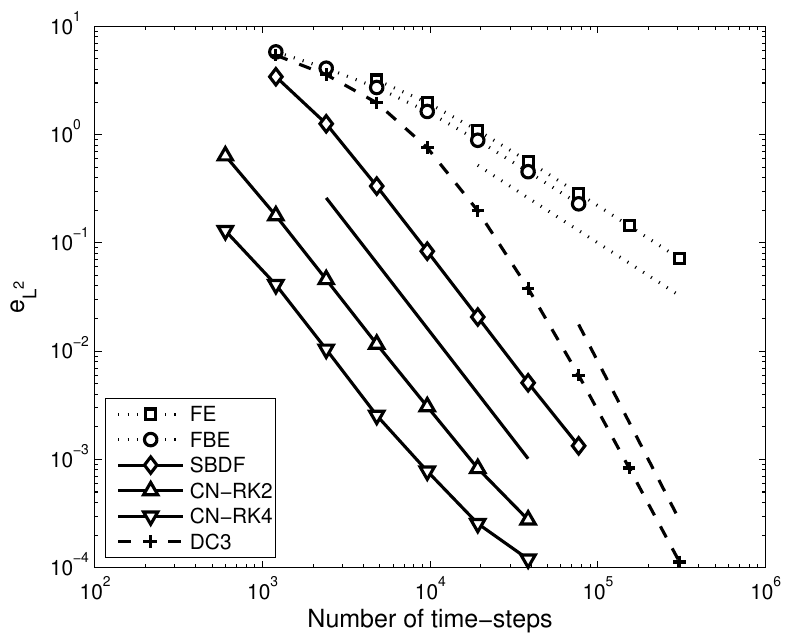}

  \caption{Error in $L^2(\Omega)$ norm} \label{fig:MSL2}
  \end{subfigure}
   \hfill
  \begin{subfigure}[b]{0.49\textwidth}
\includegraphics[width=\linewidth]{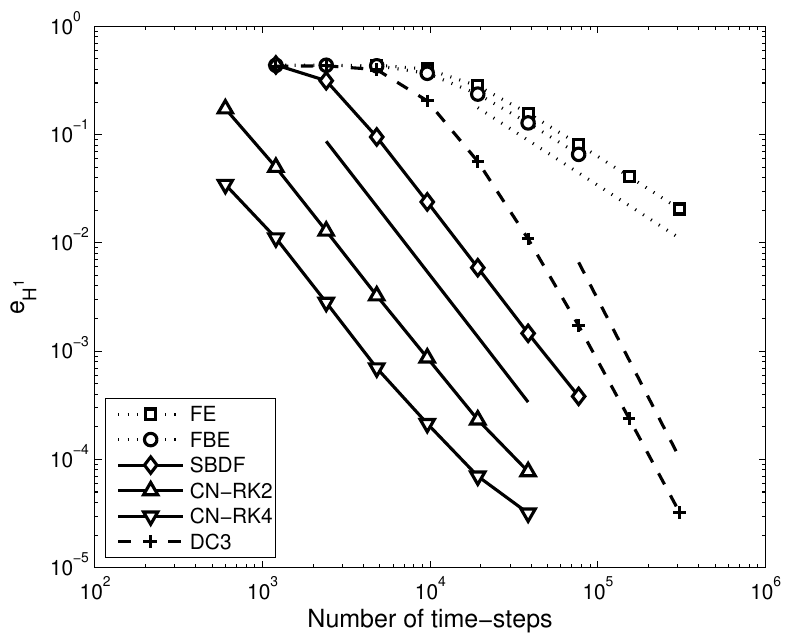}
 
  \caption{Error in $H^1(\Omega)$ seminorm} \label{fig:MSH1}
\end{subfigure}
 \vskip\baselineskip
 \begin{subfigure}[b]{0.49\textwidth}
\includegraphics[width=\linewidth]{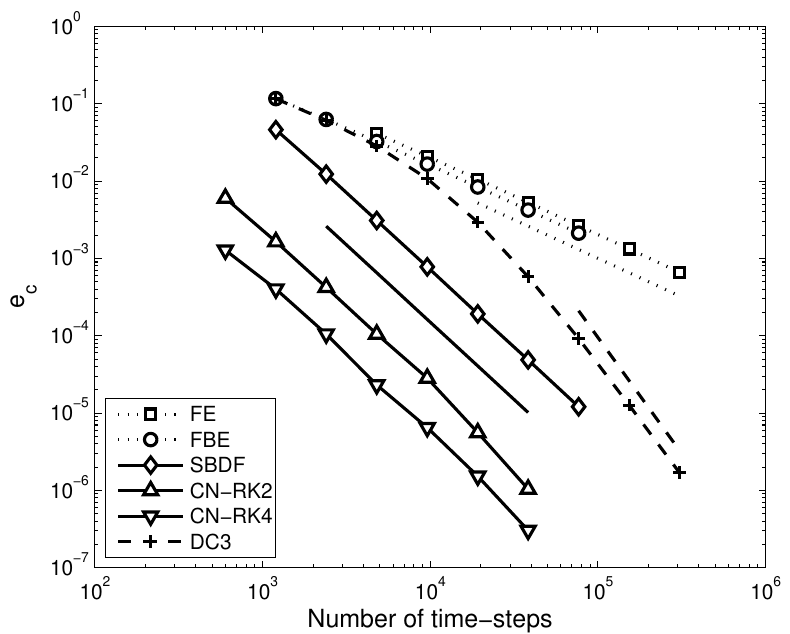}

  \caption{Error of the wave velocity $c$} \label{fig:MSc}
  \end{subfigure}
   \hfill
  \begin{subfigure}[b]{0.49\textwidth}
\includegraphics[width=\linewidth]{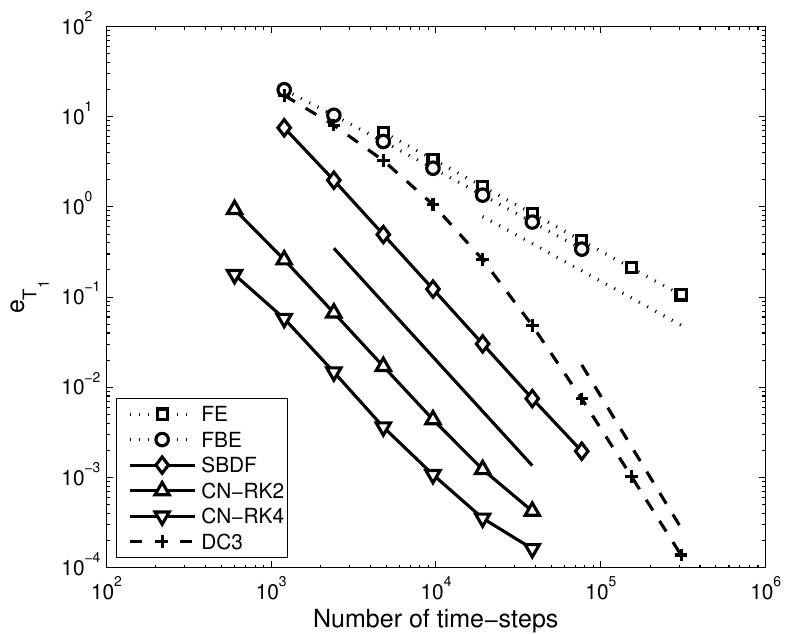}
 
  \caption{Error of the depolarization time $T_1$} \label{fig:MST}
\end{subfigure}

\caption{Errors for the MS model at fixed $h=1$. The order of convergence can be observed by comparing with lines of negative slope one (dotted), two (solid), and three (dashed).}\label{fig:MSconv}
  
\end{figure}

To determine $c$ and $T_1$, we take $\hat{u} = 0.5$, $x_1 = 50$ and $x_2 = 80$. For the convergence of the DC3 method, we used the same method instead of SBDF2 to calculate a reference solution with $\Delta t =7/60000$. The results of the convergence tests for the 1D MS model are shown in Figure \ref{fig:MSconv}.

For the MS model, all methods studied showed their expected asymptotic order of convergence, but the Strang splitting methods lose their order when $\Delta t$ gets smaller. This could be caused by the discontinuities in the ionic model. For the first-order methods, we observe that for the same $\Delta t$, FBE is slightly more accurate than its fully explicit version, FE.

For second-order methods, we observe that for the same $\Delta t$, CN-RK4 is almost ten times more accurate than CN-RK2, which in turn is about twenty times more accurate than SBDF2. For the DC3 method, it takes very small values of $\Delta t$ before the method exhibits its correct order of convergence. Even for $\Delta t = 7/6144$, its error is still larger than what we would obtain by using a good second-order method.

\subsubsection{ten Tuscher-Noble-Noble-Panfilov with 1D Monodomain}
Due to the complexity and stability requirements of the TNNP model, we study the convergence of the methods on a smaller spatial domain. The domain will be too small for the whole wave to develop, but still includes the depolarization wavefront. The reference solution is computed using SBDF2 on a domain of length $\SI{5}{\centi\metre}$ discretized in space with $160$ nodes and computed at time $\SI{12}{\milli\second} $ with $\Delta t = 6e-7$. The reference solution is shown in Figure \ref{fig:refsoltnnp}. 

\begin{figure}[!htb]
\centering
\begin{subfigure}[b]{0.49\textwidth}
\includegraphics[width=\linewidth]{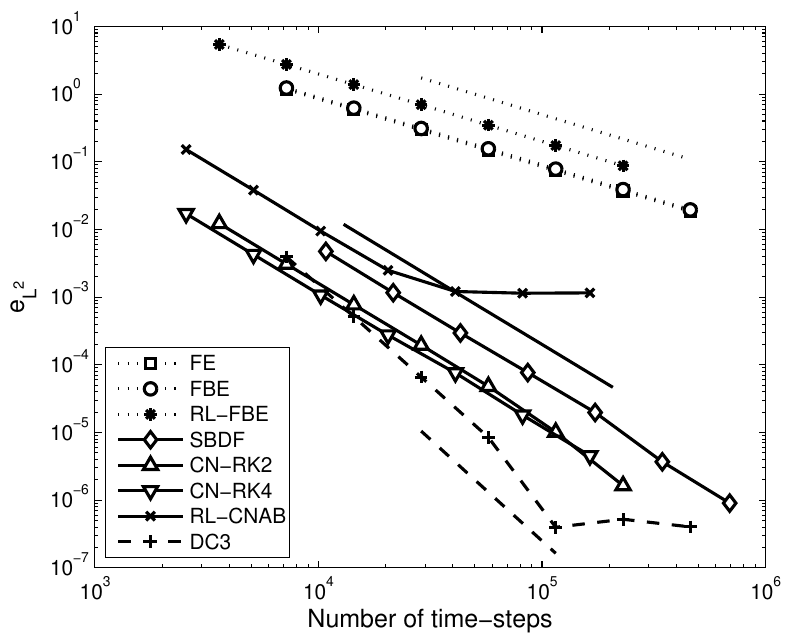}

  \caption{Error in $L^2(\Omega)$ norm} \label{fig:TNNPL2}
  \end{subfigure}
   \hfill
  \begin{subfigure}[b]{0.49\textwidth}
\includegraphics[width=\linewidth]{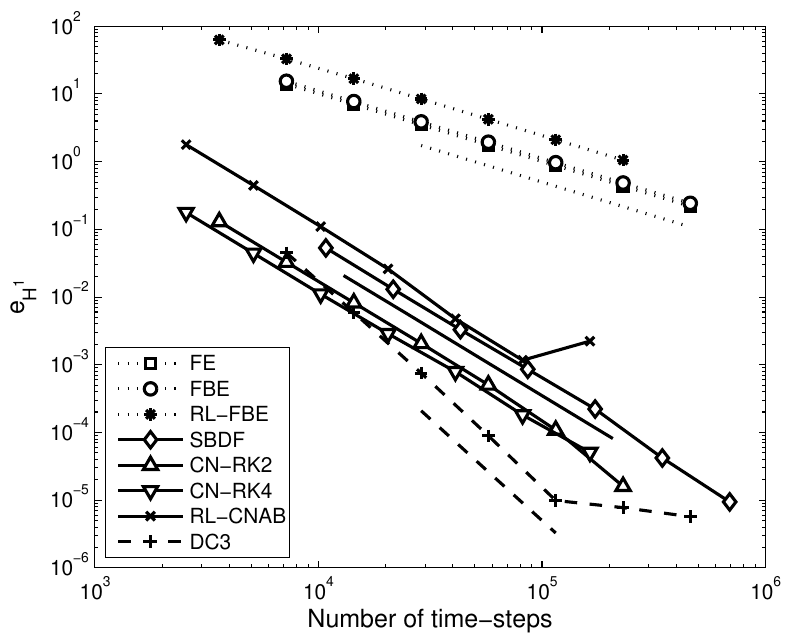}
 
  \caption{Error in $H^1(\Omega)$ seminorm} \label{fig:TNNPH1}
\end{subfigure}
 \vskip\baselineskip
 \begin{subfigure}[b]{0.49\textwidth}
\includegraphics[width=\linewidth]{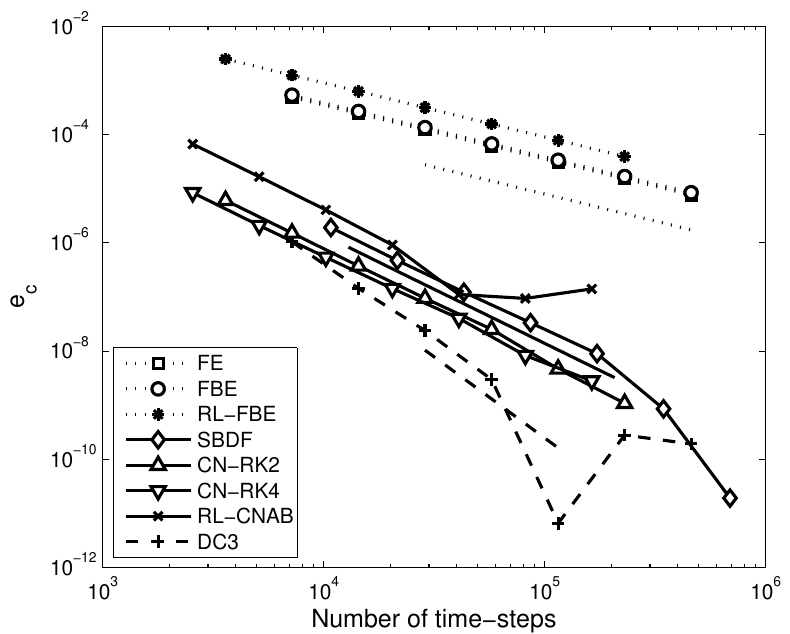}

  \caption{Error of the wave velocity $c$} \label{fig:TNNPc}
  \end{subfigure}
   \hfill
  \begin{subfigure}[b]{0.49\textwidth}
\includegraphics[width=\linewidth]{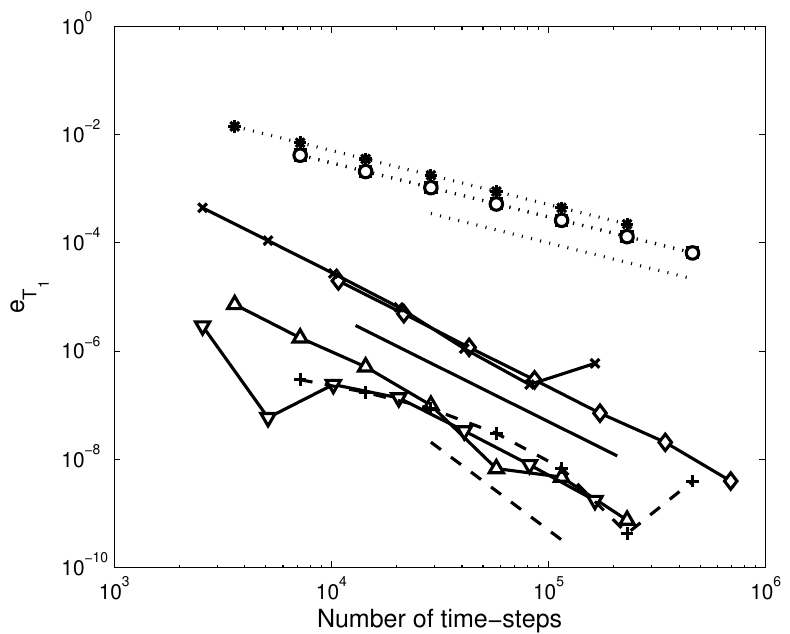}
 
  \caption{Error of the depolarization time $T_1$} \label{fig:TNNPT}
\end{subfigure}

\caption{Errors for the TNNP model at fixed $h=1/32$. The order of convergence can be observed by comparing with lines of negative slope one (dotted), two (solid), and three (dashed).}\label{fig:TNNPconv}
  
\end{figure}

To determine $c$ and $T_1$, we take $\hat{u} = \SI{-30}{\milli \volt}$, $x_1 = \SI{1}{\centi\meter}$ and $x_2 = \SI{2.5}{\centi\meter}$. For the convergence of the DC3 method, we used the same method instead of SBDF2 to calculate a reference solution with $\Delta t =7.5e-7$. The results of the convergence tests for the 1D TNNP model are shown in Figure \ref{fig:TNNPconv}.

For the TNNP model, all methods studied showed their expected asymptotic order of convergence, but the convergence of the RL-CNAB and DC3 methods is erratic when $\Delta t$ gets very small. The solution obtained with the RL-CNAB method first converges to the reference solution obtained with SBDF2 then the convergence stagnates when $\Delta t$ gets very small. The $L^2$ error between the solutions seems to stabilize around $e_{L^2} = 0.00115$, which is less than $0.0012\%$ of the reference solution's $L^2$ norm. The DC3 method shows similar trends with third-order convergence that deteriorates when $\Delta t$ gets very small. This is likely due to the difficulty of computing a reference solution at such a small level of errors. Another reason for this strange behaviour with very small time-steps might be caused by the discontinuities in the right-hand-side of the ODEs for the TNNP model. For a small time-step, it is more likely that the numerical solution falls close to a discontinuity at some of the grid points. Before its order of convergence deteriorates, the DC3 method is actually more accurate than the second-order methods.

It is important to note that the TNNP model is very stiff and that is why we use very small time-steps for the convergence tests. In practice, it is not relevant to have errors as small as those obtained for the smallest time-steps used since the modelling error is then much larger then the numerical error.

For the first-order methods, we observe that for the same $\Delta t$, FE is slightly more accurate than its semi-implicit version, FBE, and both methods are more than twice as accurate as RL-FBE. 

For second-order methods, the $\Delta t$ used are not the same, but one can easily see that for the same $\Delta t$, CN-RK4 is almost two times more accurate than CN-RK2, which in turn is about two times more accurate than SBDF2. This last method is about two times more accurate than RL-CNAB.

\subsubsection{Beeler-Reuter with 2D Monodomain}
\begin{figure}[htb!]
 \centerline{\includegraphics[scale=0.20]{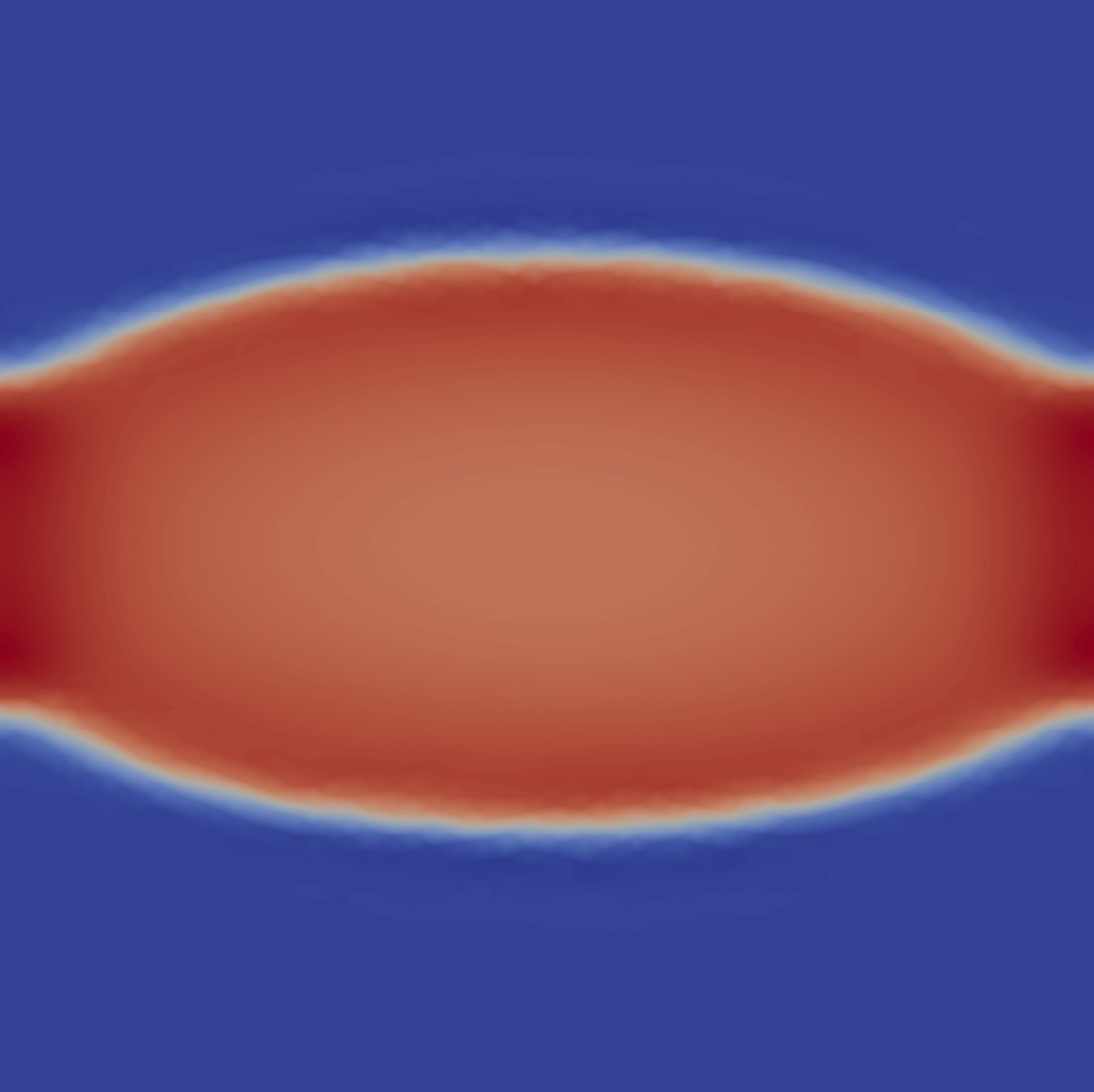}}
 
  \caption{Reference solution for the transmembrane potential at time $T=\SI{16}{\milli\second}$ for the BR model in 2D plotted as a function of $x\in [0,1]\times [0,1]$.} \label{fig:refsol2D}
\end{figure}

In the 2D case, the reference solution is computed using SBDF2 on a $1\si{\centi\metre}\times 1\si{\centi\metre}$ square domain discretized with an unstructured mesh of 3432 points (roughly of size $59\times 59$) and computed at time $T = \SI{16}{\milli\second} $ with $\Delta t = 5.98e-6$. The reference solution is shown in Figure \ref{fig:refsol2D}.

\begin{figure}[!htb]
\centering
\begin{subfigure}[b]{0.49\textwidth}
\includegraphics[width=\linewidth]{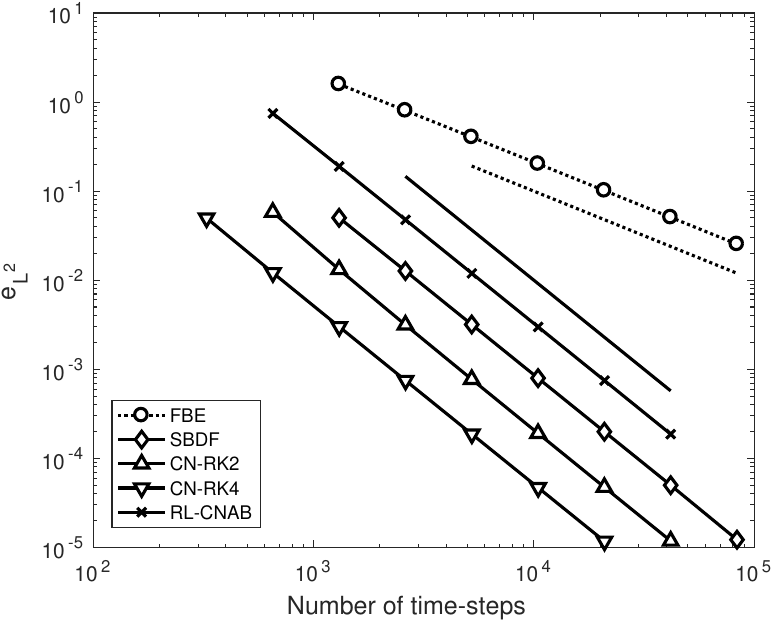}

  \caption{Error in $L^2(\Omega)$ norm} \label{fig:BR2DL2}
  \end{subfigure}
   \hfill
  \begin{subfigure}[b]{0.49\textwidth}
\includegraphics[width=\linewidth]{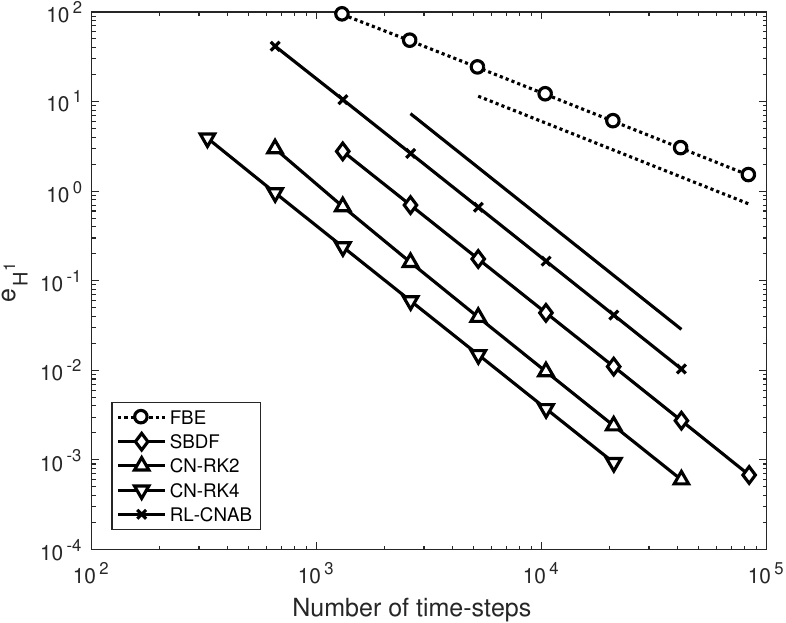}
 
  \caption{Error in $H^1(\Omega)$ seminorm} \label{fig:BR2DH1}
\end{subfigure}

\caption{Errors for the BR model in 2D. The order of convergence can be observed by comparing with lines of negative slope one (dotted), and two (solid).}\label{fig:BR2Dconv}
  
\end{figure}

The results of the convergence tests for the 2D BR model are shown in Figure \ref{fig:BR2Dconv}. We present errors and convergence rates in $L^2$ norm and $H^1$ seminorm only; wave velocity and depolarization time were not considered. We also only considered a subset of the methods studied in the 1D case.

As for the 1D BR model, all of the methods studied showed their expected asymptotic order of convergence for the 2D case. For the same $\Delta t$, CN-RK4 is about five times more accurate than CN-RK2, which in turn is almost four times more accurate than SBDF2. This last method and RL-CNAB have approximately the same accuracy. The first-order FBE method is significantly less accurate than the higher-order methods. \\

\subsection{CPU Performance of the Numerical Methods}\label{sec:cpu}

In the last sections, we compared the accuracy of the methods for given time-steps. However, this does not take into account the difference of computational cost for an iteration depending on the method used. 

We now compare the accuracy of the different numerical methods studied with respect to the computational time of the simulation. For all 1D test cases, we use the same Matlab scripts, spatial domain and spatial discretization as for the stability and convergence tests on a Lenovo ThinkCentre M900 with 4 Intel Core i5-6500T 2.10GHz processors and 7.7 GB of RAM. The 2D/3D cases were done with the compiled Fortran 2003 code CHORAL \cite{Choral} on a HP ZBook 17 with 4 Intel i7 2.70GHz processors and 16 GB of RAM. For the 2D simulations, we used the same $1\si{\centi\metre}\times 1\si{\centi\metre}$ square domain as above discretized with 3241 grid points. For the 3D simulations, we used a $1\si{\centi\metre}\times 1\si{\centi\metre}\times 1\si{\centi\metre}$ cubic domain discretized with 109620 grid points. In 2D/3D cases, the final time $T$ was set to $\SI{20}{\milli\second}$
for the BR model and $\SI{15}{\milli\second}$ for the TNNP model since the wave travels faster in this latter case. The applied stimulation current is given by (\ref{Istim2D}) with $t_0=4\si{\milli\second}$,  $\tau_0=1.5\si{\milli\second}$ and $r_0=0.125\si{\centi\metre}$.

For the sake of comparing the relative errors for different norms, we give the $L^2$ norm and $H^1$ seminorm for the different reference solutions. These are, respectively, 519.6  and 155.1 for the 1D BR model, 19.00 and 0.3154 for the 1D MS model, 96.26 and 201.8 for the 1D TNNP model, 134.3 and 728.9 for the 2D BR model, 141.5 and 876.4 for the 2D TNNP model, 223.2 and 775.2 for the 3D BR model.

We choose a target value of the $L^2$ error on the numerical solution, $e_{L^2}$, and for each method, we find the time-step needed to achieve this level of accuracy and evaluate the computational time of the simulation. For all test cases except the BR and TNNP models in 1D, we present results for two target values of the $L^2$ error, corresponding to about $0.1\%$ and $1\%$ relative $L^2$ error. In 1D, the relative errors were set to $0.5\%$ and $0.005\%$ for BR and TNNP models, respectively. The target relative error was set to a lower value for the 1D TNNP model to present at least one test case where most time-stepping methods are within their stability region. These percentages of relative of relative errors ($0.1\%$ to $1\%$) are representative of the level of numerical error usually expected from computations in cardiac electrophysiology, given  that ionic models hardly achieve this level of accuracy. These percentages also ensure that the time steps used are as close as possible to values used in common cardiac electrophysioloy simulations. When a particular method is not stable for a given percentage of relative error, results are presented for the largest time step that gives a stable solution.

Results are shown for the different methods in Tables \ref{tab:CPUBR} to \ref{tab:CPUTNNP2D}. Each table contains the $L^2$ error, $e_{L^2}$, close enough to the chosen target error, with its corresponding $H^1$ error, $e_{H^1}$, the time-step $\Delta t$ required, the total CPU time of the simulation and the average CPU time per iteration. We indicate in the title of each table the relative size of the error $e_{L^2}$ with respect to the $L^2$ norm of the reference solution. Note that the relative $H^1$ seminorm errors usually differ from the target percentage of relative $L^2$ error, and are often greater than their respective relative $L^2$ errors.

\begin{table}[htb!]
\caption{CPU time of the numerical methods for the BR model in 1D for $0.5\%$ relative $L^2$ error}\label{tab:CPUBR}
\sisetup{round-mode=figures, round-precision=4,table-parse-only}
\centering
\begin{tabular}{|l|S|S|S|S|S|}
\hline
Method & {$e_{L^2}$} & {$e_{H^1}$} &  {$\Delta t$} & {CPU time (\si{\second})} & {CPU/time-step (\si{\milli \second})} \\ \hline\hline

 SBDF2  & 2.59048923991387 & 18.6731772114248   & 0.01         & 37.213122      & 0.93032805   \\ \hline
 RL-CNAB& 2.60384056438492 & 18.7422267907901   & 0.0053333333 & 86.097275      & 1.14796366 \\ \hline
 CN-RK4 & 2.59371736142417 & 18.5321998635816   & 0.0454545455 & 57.105731      & 6.48928761    \\ \hline
 CN-RK2 & 2.59094703842188 & 18.6291192879054   & 0.016        & 77.009233      & 3.08036932 \\ \hline
 FBE    & 25.971261229217  & 153.014626378104   & 0.0026936027 & 117.785241     & 0.79316666 \\ \hline
 DC3    & 2.605761141561   & 18.763830817885    & 0.0014925373 & 780.163321     & 2.91105716 \\ \hline
\end{tabular}
\end{table}

For the 1D BR model, we see that the most efficient method is the SBDF2 method. It is almost twice as fast as the next most efficient method, CN-RK4, which is itself faster than the CN-RK2 method. For each time-step of the SBDF2 method, we only have to compute the functions of the ionic model once, as opposed to the Strang-splitting methods where we compute them four or eight times for CN-RK2 and CN-RK4, respectively. As for the DC3 method, the functions are computed only three times per time-step, but it takes small values of $\Delta t$ for the method to reach a higher rate of convergence. We observe that the third-order method, DC3, is not efficient compared to the second-order methods. Indeed, it takes more than twenty times more CPU time to compute the solution than the SBDF2 method for the chosen level of error. By choosing very small errors, we expect the DC3 method to eventually surpass the second-order methods. At the chosen level of error $e_{L^2}$, the DC3 method has not yet entered its asymptotic zone. 

We only tested the computational time of the FBE method at $5 \%$ relative $L^2$ error because it would take extremely small time-steps to reach the level of error that we chose for the other methods. It is clear that for $0.5\%$ relative error, the FBE method is not efficient compared to high-order methods. 

The computational time per time-step of the SBDF2 method is only slightly larger than for the FBE method and thus there is no advantage in terms of efficiency to use a first-order method. One time-step of CN-RK4 is two times more costly than for CN-RK2, which is almost four times more costly than the SBDF2 method. This is directly related to the number of computations of the ionic functions. One time-step of DC3 is more than three times the cost of an SBDF2 time-step, which also relates to the computation of the ionic functions, but there is additional cost due to the more complex nature of the DC3 method. These relations extend to the different models studied. Even though RL-CNAB also only has one computation of the ionic functions per time-step, the use of the exponential function makes an iteration for this method slightly more expensive than for the SBDF2 method. This impact is lessened in the case of the TNNP model because the computation of the ionic functions takes up most of the computational cost, due to the complexity of the TNNP model.

\begin{table}[htb!]
\caption{CPU time of the numerical methods for the MS model in 1D for $0.1\%$ $L^2$ error}
\centering
\begin{tabular}{|l|c|c|c|c|c|}
\hline
Method  & {$e_{L^2}$} & {$e_{H^1}$} &  {$\Delta t$} & {CPU time (\si{\second})}  & {CPU/time-step (\si{\milli \second})} \\ \hline\hline

SBDF2  & 0.01898 & 0.005431 & 0.01728  & 2.670 & 0.1319\\ \hline
CN-RK4 & 0.01884 & 0.005091 & 0.1971   & 1.597 & 0.8990\\ \hline
CN-RK2 & 0.01894 & 0.005329 & 0.09272  & 1.736 & 0.4599\\ \hline
FBE    & 0.1897  & 0.05432  & 0.003763 & 12.09 & 0.1300\\ \hline
DC3    & 0.01903 & 0.005461 & 0.006972 & 23.25 & 0.4671\\ \hline

\end{tabular}
\end{table}

\begin{table}[htb!]
\caption{CPU time of the numerical methods for the MS model in 1D for $1\%$ $L^2$ error}
\sisetup{round-mode=figures, round-precision=4,table-parse-only}
\centering
\begin{tabular}{|l|c|c|c|c|c|}
\hline
Method  & {$e_{L^2}$} & {$e_{H^1}$} &  {$\Delta t$} & {CPU time (\si{\second})}  & {CPU/time-step (\si{\milli \second})} \\ \hline\hline

 SBDF2       & 0.18997 & 0.05426 & 0.05488  & 0.9671 & 0.1516  \\ \hline
 CN-RK4      & 0.18970 & 0.05093 & 0.7447   & 0.4793 & 1.019   \\ \hline
 CN-RK2      & 0.18855 & 0.05276 & 0.3017   & 0.6061 & 0.5225  \\ \hline
 FBE         & 0.18991 & 0.05437 & 0.003767 & 13.75  & 0.1480  \\ \hline
 DC3         & 0.18991 & 0.05422 & 0.01766  & 10.22  & 0.5157  \\ \hline

\end{tabular}
\end{table}

For the 1D MS model, we see that the most efficient method is the CN-RK4 method, followed closely by the CN-RK2 method. The CN-RK4 method is almost twice as efficient as the SBDF2 method. We observe that our third-order method, DC3, is not efficient compared to the second-order methods. Indeed, it takes almost 8-10 times more CPU time to compute the solution than the SBDF2 method for the chosen error. These conclusions hold for both levels of relative error.

We only tested the computational time of the FBE method at $1 \%$ relative $L^2$ error because it would take extremely small time-steps to reach the $0.1\%$ relative error level. It is easy to figure that for $0.1\%$ relative error, the FBE method is not efficient compared to high-order methods.

\begin{table}[htb!]
\caption{CPU time of the numerical methods for the TNNP model in 1D for $0.005\%$ relative $L^2$ error}
\sisetup{round-mode=figures, round-precision=4,table-parse-only}
\centering
\begin{tabular}{|l|S|S|S|S|S|}
\hline
Method & {$e_{L^2}$} & {$e_{H^1}$} &  {$\Delta t$} & {CPU time (\si{\second})} & {CPU/time-step (\si{\milli \second})}  \\ \hline\hline

 SBDF2  & 0.00480000893782859 & 0.0538635709797139 & 0.0011162791 & 4.249266 & 0.395280558    \\ \hline
 RL-CNAB& 0.00480189858210097 & 0.0546684640000704 & 0.0008344924 & 6.862648   & 0.4772356    \\ \hline
 CN-RK4 & 0.00478544030906128 & 0.0495178366845316 & 0.0024870466 & 14.334048  & 2.9707871 \\ \hline
 CN-RK2 & 0.00479746280420562 & 0.0508413756514704 & 0.0020779221 & 12.5424804000013  & 2.1718580779 \\ \hline
 FBE    & 0.0479627158696248  & 0.596785139907064  & 6.41711229946524E-005      & 73.473536             & 0.3929066\\ \hline
 DC3    & 0.003961010199763   & 0.045096108751937  & 0.00166666             & 9.388821  & 1.30400 \\ \hline

\end{tabular}
\end{table}

For the 1D TNNP model, most methods require very small time-steps to have stable solutions. Therefore the relative error is very small compared to those of the previous models. We see that the most efficient method is, as for BR, the SBDF2 method. For $0.005\%$ relative $L^2$ error, it is about 1.5 times faster than the next most efficient second-order method, RL-CNAB, nearly three times as fast as the CN-RK2 method, and 3.4 times faster than the CN-RK4 method. Interestingly, we observe that the DC3 method is the third most efficient method, since DC3 is unexpectedly accurate for TNNP.

We only tested the computational time of the FBE method at $0.05 \%$ relative $L^2$ error because it would take extremely small time steps to reach the error we chose for the other methods. Because of the stability constraints on the other methods, only the RL-CNAB can have an error as high as $0.05\%$. It would be easy to show that when RL-CNAB is in its asymptotic zone, it is more efficient than FBE. It is also easy to figure that for $0.005\%$ relative error, the FBE method is not efficient compared to high-order methods.

\begin{table}[htb!]
\caption{CPU time of the numerical methods for the BR model in 2D for $0.1\%$  $L^2$ error}
\sisetup{round-mode=figures, round-precision=4,table-parse-only}
\centering
\begin{tabular}{|l|S|S|S|S|S|}
\hline
Method & {$e_{L^2}$} & {$e_{H^1}$} &  {$\Delta t$} & {CPU time (\si{\second})}  & {CPU/time-step (\si{\milli \second})} \\ \hline\hline

SBDF2	& 0.129067391 	& 5.08239222	& 0.0108108111	& 1.62047696	& 0.8759334918\\ \hline
CN-RK4	& 0.128501192	& 4.22034121	& 0.0235294122	& 3.24099517	& 3.8129354941\\ \hline
CN-RK2	& 0.131661892	& 4.48622799	& 0.0178571437	& 2.65156102	& 2.3674651964\\ \hline
DC2	& 0.133484393	& 4.98323965	& 0.00424628472	& 8.42094803	& 1.7878870552\\ \hline
DC3	& 0.134858370	& 5.10799789	& 0.0108108111	& 5.10576582	& 2.7598734162\\ \hline
RL-CNAB	& 0.132911608 	& 5.05598402	& 0.00677966094	& 3.03639603 	& 1.0292867898\\ \hline
FBE 	& 0.133849546	& 5.08005762	& 0.000500000024	& 30.1190720	& 0.7529768\\ \hline
RL-FBE 	& 0.134287104	& 5.01957607	& 0.000173913038	& 98.5357361	& 0.8568324878\\ \hline

\end{tabular}
\end{table}

\begin{table}[htb!]
\caption{CPU time of the numerical methods for the BR model in 2D for $1\%$  $L^2$ error}\label{tab:CPUBR2D}
\sisetup{round-mode=figures, round-precision=4,table-parse-only}
\centering
\begin{tabular}{|l|S|S|S|S|S|}
\hline
Method & {$e_{L^2}$} & {$e_{H^1}$} &  {$\Delta t$} & {CPU time (\si{\second})}  & {CPU/time-step (\si{\milli \second})} \\ \hline\hline

SBDF2	& 0.264444798 	& 10.3999052	& 0.0156250000	& 1.04476190	& 0.816220234375\\ \hline
CN-RK4	& 0.922152758	& 30.2841969	& 0.0625000000	& 1.19145393	& 3.72329353125\\ \hline
CN-RK2	& 0.708067656	& 24.0508232	& 0.0399999991	& 1.12998390	& 2.2599678\\ \hline
DC2	& 1.27236688	& 46.9494781	& 0.0156250000	& 2.02250910	& 1.580085234375\\ \hline
DC3	& 0.579459250	& 21.7128201	& 0.0199999996	& 2.44924498	& 2.44924498\\ \hline
RL-CNAB & 1.12534392	& 42.3262672	& 0.0199999996	& 0.951569080 	& 0.951569080\\ \hline
FBE 	& 1.34418201	& 50.4715347	& 0.00499999989	& 2.97760105	& 0.7444002625\\ \hline
RL-FBE  & 1.23079884	& 45.4310417	& 0.00159999996	& 10.5121183 	& 0.840969464\\ \hline

\end{tabular}
\end{table}

\begin{table}[htb!]
\caption{CPU time of the numerical methods for the BR model in 3D for $0.1\%$  $L^2$ error}
\sisetup{round-mode=figures, round-precision=4,table-parse-only}
\centering
\begin{tabular}{|l|S|S|S|S|S|}
\hline
Method & {$e_{L^2}$} & {$e_{H^1}$} &  {$\Delta t$} & {CPU time (\si{\second})}  & {CPU/time-step (\si{\milli \second})} \\ \hline\hline

SBDF2	& 0.113516435 	& 3.77564859	& 0.00999999978	& 137.387543	& 68.6937715\\ \hline
CN-RK4	& 0.221395016	& 7.18596458	& 0.0245398767	& 168.453735	& 206.6916993\\ \hline
CN-RK2	& 0.222662464	& 7.25207376	& 0.0186915882	& 167.500275	& 156.54231308\\ \hline
DC2	& 0.220239803	& 7.28450584	& 0.00438596494	& 703.194580	& 154.2093377\\ \hline
DC3	& 0.227338701	& 7.51837730	& 0.0116618080	& 436.441772	& 254.48499825\\ \hline
RL-CNAB	& 0.222334906	& 7.38102341	& 0.00760456268	& 192.265701 	& 73.104829277\\ \hline
FBE 	& 0.259627044	& 8.92350388	& 0.000475059380& 2089.74658	& 49.6376859857\\ \hline
RL-FBE 	& 0.277068824	& 9.31719017	& 0.000177541049& 5555.35107	& 49.31514487\\ \hline

\end{tabular}
\end{table}

\begin{table}[htb!]
\caption{CPU time of the numerical methods for the BR model in 3D for $1\%$  $L^2$ error}
\sisetup{round-mode=figures, round-precision=4,table-parse-only}
\centering
\begin{tabular}{|l|S|S|S|S|S|}
\hline
Method & {$e_{L^2}$} & {$e_{H^1}$} &  {$\Delta t$} & {CPU time (\si{\second})}  & {CPU/time-step (\si{\milli \second})} \\ \hline\hline

SBDF2	& 0.113520928 	& 3.77585626	& 0.00999999978	& 139.009811	& 69.5049055\\ \hline
CN-RK4	& 0.360116541	& 11.6975975	& 0.0312500000	& 133.819641	& 209.0931890625\\ \hline
CN-RK2	& 0.643448114	& 20.9806576	& 0.0312500000	& 102.984673	& 160.9135515625\\ \hline
DC2	& 1.95682395	& 64.1476898	& 0.0156250000	& 198.716492	& 155.247259375\\ \hline
DC3	& 0.463779598	& 15.3078480	& 0.0156250000	& 318.175232	& 248.5744\\ \hline
RL-CNAB	& 2.25324059	& 74.7507782	& 0.0243902430	& 63.2703781	& 77.15899768\\ \hline
FBE 	& 2.24631286	& 76.9409256	& 0.00480769249	& 257.744568	& 61.957828846\\ \hline
RL-FBE 	& 2.23764157	& 75.0969086	& 0.00178571430	& 676.477783	& 60.39980205\\ \hline

\end{tabular}
\end{table}

For the 2D case with the BR model, all methods were stable enough to achieve a relative error of $0.1\%$. At this level of error, the SBDF2 method is, as in the 1D case, the most competitive, followed by RL-CNAB and the CN-RKn methods. The DCn methods are not competitive
as they are not yet in their asymptotic region of convergence. The first order methods FBE and RL-FBE are not recommended. For a relative error of $1\%$, only the FBE, DC2 and RL methods can provide solutions at this exact error level. All other methods had their time step restricted by stability, resulting in error levels  below $1\%$. Nevertheless, the SBDF2 and CN-RKn methods remains competitive with RL-CNAB, which is marginally more efficient in terms of CPU time. Accounting for the fact that SBDF2 is about four times more accurate for about the same CPU time as RL-CNAB, this is again the winning method. The RL-FBE method is to be avoided due to its lack of accuracy, requiring a time step ten times smaller than SBDF2 to reach $1\%$ error. In 3D, both with a relative error of $0.1\%$ and $1\%$, all the comments made above for the 2D case on stability, accuracy and efficiency ordering of the methods remain valid, except that at $1\%$ error the RL-CNAB method is twice as fast as SBDF2 but with a $L^2$ error twenty times larger. Also, at $1\%$ error, the gap in efficiency between first and second order methods is reduced.

\begin{table}[htb!]
\caption{CPU time of the numerical methods for the TNNP model in 2D for $0.1\%$  $L^2$ error}
\sisetup{round-mode=figures, round-precision=4,table-parse-only}
\centering
\begin{tabular}{|l|S|S|S|S|S|}
\hline
Method & {$e_{L^2}$} & {$e_{H^1}$} &  {$\Delta t$} & {CPU time (\si{\second})}  & {CPU/time-step (\si{\milli \second})} \\ \hline\hline

SBDF2	& 0.00216718810 & 0.0992370397	& 0.000600000028	& 40.5011787	& 1.620047148\\ \hline
CN-RK4	& 0.00358875026	& 0.156160071	& 0.00200000009		& 58.1130791	& 7.748410667\\ \hline
CN-RK2	& 0.00318990485	& 0.140044570	& 0.00156250002		& 42.7543335	& 4.45357640625\\ \hline
DC2	& 0.0197716877	& 0.894100070	& 0.000899820065	& 54.0457039	& 3.24209391975\\ \hline
DC3	& 0.000772150583& 0.0354534835	& 0.000899820065	& 84.5447464	& 5.0716706134\\ \hline
RL-CNAB	& 0.147419095	& 6.83233738	& 0.00524475519		& 5.66641235 	& 1.981263038756\\ \hline
FBE 	& 0.136346474	& 6.51243639	& 0.000403225800	& 57.9298973	& 1.5572552788\\ \hline
RL-FBE 	& 0.137163579	& 6.33085108	& 0.000154400419	& 169.245789	& 1.742108049039\\ \hline

\end{tabular}
\end{table}

\begin{table}[htb!]
\caption{CPU time of the numerical methods for the TNNP model in 2D for $1\%$  $L^2$ error}\label{tab:CPUTNNP2D}
\sisetup{round-mode=figures, round-precision=4,table-parse-only}
\centering
\begin{tabular}{|l|S|S|S|S|S|}
\hline
Method & {$e_{L^2}$} & {$e_{H^1}$} &  {$\Delta t$} & {CPU time (\si{\second})}  & {CPU/time-step (\si{\milli \second})} \\ \hline\hline

SBDF2	& 0.00216718810 & 0.0992370397	& 0.000600000028	& 40.5011787	& 1.620047148\\ \hline
CN-RK4	& 0.00358875026	& 0.156160071	& 0.00200000009		& 58.1130791	& 7.748410667\\ \hline
CN-RK2	& 0.00318990485	& 0.140044570	& 0.00156250002		& 42.7543335	& 4.45357640625\\ \hline
DC2	& 0.0197716877	& 0.894100070	& 0.000899820065	& 54.0457039	& 3.24209391975\\ \hline
DC3	& 0.000772150583& 0.0354534835	& 0.000899820065	& 84.5447464	& 5.0716706134\\ \hline
RL-CNAB	& 1.42418075	& 65.8968735	& 0.0174418613		& 1.82400799 	& 2.12093952325\\ \hline
FBE 	& 0.317001373	& 15.1256075	& 0.000949367066	& 25.1384048	& 1.59103827848\\ \hline
RL-FBE 	& 1.32150400	& 60.1157494	& 0.00156250002		& 18.2760525	& 1.90375546875\\ \hline

\end{tabular}
\end{table}

For the 2D case with the TNNP model, the only methods that could provide a solution with a relative error of $0.1\%$ were FBE and RL methods, while at $1\%$ error only the RL methods are stable enough. FBE can provide solutions with such error levels not because it is more stable. FBE simply requires such a small time step for accuracy that this one falls within its stability region.
The clear winner for CPU time is the RL-CNAB method. Among the other higher order methods, SBDF2 wins in terms of CPU time, this one being seven times larger than for RL-CNAB at $0.1\%$ error but with a $L^2$ error 70 times smaller than for RL-CNAB. Comparing our 1D and 2D results for the TNNP model, second order methods are nearly equally competitive at low error levels, but it is obvious that the lack of stability of non-RL methods makes them less attractive for solving the stiffest ionic models with larger numerical error levels and affordable time steps.

The difference in the performance of the methods between the 1D and 2D/3D cases is due to the relative cost of solving the linear systems in 1D versus 2D/3D, compared to the cost of evaluating the ionic functions. Besides the different linear solvers used, the implementation is done with MATLAB in 1D versus the compiled language Fortran for the 2D/3D code. In 1D, we use a direct solver with a tridiagonal matrix that makes the cost of solving the linear system about $1/16$ of doing one evaluation of the BR reaction terms. In 2D, we use an iterative solver that makes the cost of solving the linear system about 3 times more expensive than doing one evaluation of the reaction terms.


\section{Conclusion}
In this article, we provide a strategy to practically determine precise values of the critical time step for linearly implicit time-stepping methods applied to reaction-diffusion equations. The methodology requires solving numerically the reaction-diffusion equation in dimension one, but gives critical time steps that equally applies to n-dimensional settings. The efficiency of the methodology is illustrated using the monodomain model with ionic models of varying stiffness.

We also investigated the impact of ionic model complexity and stiffness on finding an efficient numerical solution through a variety of time-stepping methods. 
The simplest model that we tested is the MS model. Its stiffness is low and hence there is no need for very stable methods to be used. The most accurate method to solve this problem is Strang splitting using Runge-Kutta methods (RK2 or RK4) to solve the ionic model and the reaction part of the monodomain model, and Crank-Nicholson (CN) for the diffusion part of the monodomain model. Due to the nature of operator splitting methods, they become less competitive when solving more complex models. Indeed, for each time-step, we do two RK substeps which, for the second-order RK method, implies computing the ionic currents four times, and for the original fourth order RK, eight times. As the models get more complex and realistic, computing the ionic currents takes up most of the computational cost of the algorithms. Therefore, multistep methods such as SBDF2 become more efficient with respect to computational time when solving more complex models because the ionic currents are only computed once per time-step. For the Beeler-Reuter model, the SBDF2, RL-CNAB and the Strang splitting methods are the most efficient unless a coarse solution is needed. The easier implementation of the SBDF2 scheme makes it faster with interpreted languages such as MATLAB. For the TNNP model, the SBDF2 method outclasses the other methods in terms of efficiency when a very accurate solution is required, while the RL-CNAB methods clearly wins for more affordable simulations at larger error levels. Therefore, the choice of the time-stepping method depends on the stiffness of the model. Strang splitting methods are best for the models with less stiffness, while SBDF2 becomes better as the stiffness increases and RL-CNAB is the obvious choice for the stiffest models where stability is critical. The third-order DC3 method is the most efficient only in the case of extremely accurate solutions (not shown here). In practice, simulations do not require this much accuracy as the balance between modelling and numerical error allows for a lower level of numerical accuracy.

For the models studied here, using the largest possible time-steps needed for the stability of the semi-implicit methods results in solutions that are accurate enough relative to modelling errors. If this level of accuracy is satisfactory and the remaining goal is to reduce computational time, the stability of the methods is thus the important factor. The stiffness of the more complex models calls for more stable methods, such as Rush-Larsen methods. As we have seen in Section \ref{chap:stab}, these are almost unconditionally stable. However, in Section \ref{chap:convergence}, we saw that they are less accurate than other semi-implicit methods such as SBDF2 or Strang splitting methods. For the stiffer models such as TNNP, a lower level of accuracy is not available with the less stable SBDF2 and Strang splitting methods, because even at the largest possible time-step necessary for stability, the errors are extremely low. Such accuracy is not necessary in most cases, for which Rush-Larsen type methods with larger time-steps are an obvious winner in terms of computational time.

\bibliographystyle{plain}
\bibliography{biblio}

\end{document}